\documentclass[11pt]{article}
\usepackage{mathrsfs}
\usepackage{color}
\usepackage{amsfonts}
\usepackage{amssymb}
\usepackage{bm}
\usepackage{amsmath}

\def\sqr#1#2{{\vcenter{\vbox{\hrule height.#2pt
              \hbox{\vrule width.#2pt height#1pt \kern#1pt \vrule
width.#2pt}
              \hrule height.#2pt}}}}
\def\signed #1{{\unskip\nobreak\hfil\penalty50
              \hskip2em\hbox{}\nobreak\hfil#1
              \parfillskip=0pt \finalhyphendemerits=0 \par}}
\def\endpf{\signed {$\sqr69$}}
\def\l{\lambda}

\def\t{\tau}

\def\o{\omega}
\def\p{\phi}

\def\3n{\negthinspace \negthinspace \negthinspace }
\def\2n{\negthinspace \negthinspace }
\def\1n{\negthinspace }

\def\O{\Omega}
\def\no{\noindent}

\def\ms{\medskip}

\def\q{\quad}

\def\ou{\overline{u}}
\def\ov{\overline{\mathbf{v}}}
\def\na{\nabla}
\def\ow{\overline{w}}
\def\esssup{\mathop{\rm esssup}}

\def\supp{\hbox{\rm supp$\,$}}

\def\({\Big (}
\def\){\Big )}
\def\[{\Big[}
\def\]{\Big]}

\allowdisplaybreaks

\def\be{\begin{equation}}
\def\bel{\begin{equation}\label}
\def\ee{\end{equation}}
\def\bea{\begin{eqnarray}}
\def\eea{\end{eqnarray}}
\def\bt{\begin{theorem}}
\def\et{\end{theorem}}
\def\bc{\begin{corollary}}
\def\ec{\end{corollary}}
\def\bl{\begin{lemma}}
\def\el{\end{lemma}}
\def\bp{\begin{proposition}}
\def\ep{\end{proposition}}
\def\br{\begin{remark}}
\def\er{\end{remark}}
\def\ba{\begin{array}}
\def\ea{\end{array}}
\def\bd{\begin{definition}}
\def\ed{\end{definition}}
\def\bal{\begin{align}}
\def\eal{\end{align}}

\newtheorem{lemma}{Lemma}[section]
\newtheorem{remark}{Remark}[section]

\newtheorem{theorem}{Theorem}[section]
\newtheorem{corollary}{Corollary}[section]

\newtheorem{definition}{Definition}[section]
\newtheorem{proposition}{Proposition}[section]

\makeatletter
   
   \@addtoreset{equation}{section}
\makeatother
\allowdisplaybreaks

\begin{document}

\title{\bf Exact controllability to nonnegative trajectory for a chemotaxis system
}
\date{}
\author{Qiang Tao \thanks{School of Mathematics and Statistics,
 Shenzhen University, Shenzhen 518060, China and Shenzhen Key Laboratory of Advanced Machine Learning and Applications, Shenzhen University, Shenzhen 518060, China. E-mail address:
 taoq060@126.com.\ms} \q and  \q Muming Zhang\thanks{Corresponding author. School of Mathematics and Statistics,
 Northeast Normal University, Changchun 130024, China. E-mail address:
zhangmm352@nenu.edu.cn. \ms }}

\date{}

\maketitle

\begin{abstract}
This paper studies the controllability for a Keller-Segel type chemotaxis model with singular
sensitivity. Based on the Hopf-Cole transformation, a nonlinear parabolic system, which has first-order couplings, and the coupling coefficients are functions that depend on both time and space variables, is derived. Then, the controllability result is proved by a new global Carleman estimate for general coupled parabolic equations allowed to contain a convective term.
 Also, the global existence of nonnegative solution for the chemotaxis system is discussed. 
 \end{abstract}

\noindent{\bf Key Words. } Exact controllability, chemotaxis, nonnegative trajectory, Carleman estimate

\noindent{\bf Mathematics Subject Classification (2020).} 35K55, 93B05

\section{Introduction and the main results}

Chemotaxis is a biological process in which cells move toward a chemically more favorable environment, e.g., bacteria swim to places with high concentration of food molecules.
In order to model the interaction between  vascular endothelial cells (VEC) and vascular endothelial growth factor (VEGF),
a Keller-Segel type chemotaxis model with logarithmic sensitivity was proposed in \cite{levine3}:
\begin{equation}\label{1.1}
\left\{\begin{array}{ll}
\widetilde{u}_t=\nabla\cdot(D\nabla \widetilde{u}-\chi\widetilde{u}\nabla \ln \widetilde{c}),& (\tilde{x},\tilde{t})\in\widetilde{Q},\\
\widetilde{c}_t=\Delta \widetilde{c}-\mu\widetilde{u}\widetilde{c}, &(\tilde{x},\tilde{t})\in\widetilde{Q},
\end{array}\right.
\end{equation}
where $\widetilde{\Omega}$ is a bounded domain in $\mathbb{R}^{n}$ and  $\widetilde{Q}=\widetilde{\Omega}\times (0,\widetilde{T})$ with $\widetilde{T}>0$. 
The functions $\widetilde{u}$ and $\widetilde{c}$ denote the density of VEC and concentration of VEGF, respectively.
The parameter $D>0$ is referred to as the diffusivity of VEC. The logarithmic sensitivity $\ln \widetilde{c}$ with the constant $\chi>0$ indicates that cell chemotactic responsing to VEGF follows the Weber-Fechner's law, which has many important applications in biological modelings (see \cite{AL1987, BM1985, DLK1972, KS1971}). The positive constant $\mu$ measures the degradation rate of VEGF.
The system (\ref{1.1}), which plays a central role in illustrating the spreading of cancer cells to other tissues in cancer metastasis, also indicates that the population of VEC could aggregate over time at certain
spatial locations, since it is driven against diffusion by the concentration
gradient of VEGF at spatial locations where the chemical signal increases.
This kind of aggregation may lead to the phenomenon of finite-time blow up.

It is obvious that the logarithmic sensitivity function $\ln\widetilde{c}$ is singular at $\widetilde{c}=0$.
 In order to overcome singularity, an effective approach is to apply the Hopf-Cole transformation as follows
(see, for example, \cite{levine2}):
$$\widetilde{\mathbf{v}}=\nabla\ln \widetilde{c}=\frac{\nabla\widetilde{c}}{\widetilde{c}},$$
together with  scalings $t=\frac{\chi\mu}{D}\tilde{t},\ x=\frac{\sqrt{\chi\mu}}{D}\tilde{x},\ \widehat{\mathbf{v}}(x,t)=-\sqrt{\frac{\chi}{\mu}}\mathbf{\widetilde{v}}(\tilde{x},\tilde{t})$, $\widehat{u}(x,t)=\widetilde{u}(\tilde{x},\tilde{t})$. Then the system (\ref{1.1}) is transformed into the following form:
\begin{equation}\label{h1.2}
\left\{\begin{array}{ll}
\widehat{u}_t-\Delta \widehat{u}=\nabla\cdot(\widehat{u}\widehat{\mathbf{v}}),& (x,t)\in Q,\\
\widehat{\mathbf{v}}_t-\frac{1}{D}\Delta\widehat{\mathbf{v}}=
\nabla(-|\widehat{\mathbf{v}}|^2+\widehat{u}), &(x,t)\in Q,
\end{array}\right.
\end{equation}
where $\O=\frac{\sqrt{\chi\mu}}{D}\widetilde{\O}$ and $T=\frac{\chi\mu}{D}\widetilde{T}$. For the sake of simplicity, we take $D=1$ in what follows. Here, $Q=\O\times(0,T)$, where $\O\subset\mathbb{R}^{n}(1\leq n\leq 3)$ is a bounded domain with smooth boundary $\Gamma$. Let $\Sigma=\Gamma\times (0,T)$ and $T>0$.

To prevent the spread and metastasis of tumor cells, external intervention is essential. This urges us to study the controllability of the system (\ref{h1.2}).  The controllability problem of chemotaxis models can be viewed as finding control strategies
(such as the use of drug treatment) to make the concentration of chemical and the density of cells
tend to the given substance concentration and cell density. In this paper, we will study the exact controllability of system (\ref{h1.2}) to a nonnegative trajectory defined by system (\ref{1.5}).

Let $\o$ be a given nonempty open subset of $\O$. Denote by $\chi_\o$ the characteristic function of $\o$. We will study the following controlled chemotaxis system:
\begin{equation}\label{1.2}
\left\{\begin{array}{ll}
u_t-\Delta u=\nabla\cdot(u\mathbf{v}),& (x,t)\in Q,\\
\mathbf{v}_t-\Delta\mathbf{v}=
\nabla(-|\mathbf{v}|^2+u)+\chi_{\o}\mathbf{h}, &(x,t)\in Q,\\
u=\overline{p}, ~~\mathbf{v}=0,&(x,t)\in\Sigma,\\
(u,\mathbf{v})(x,0)=(u_0,\mathbf{v}_0)(x), & x\in\Omega,
\end{array}\right.
\end{equation}
where $(u,\mathbf{v})$ is the state, $ \mathbf{h}$ is the control funtion and $\overline{p}$ is a positive constant. Obviously, the control acts on the chemical concentration equation.


In the last decades, there are many works addressing the qualitative theory of the solutions to chemotaxis models (
see for example \cite{H2003, HW2005, JLW2013, LZ2015, rwwzz, TWW2013, WXY2016, W2018}
 and the rich references therein). However, few results are known on the controllability of chemotaxis models, we refer to \cite{CG2015, CG2017, GZ2014, Guo-Zhang}.
The local null controllability for a chemotaxis system of
parabolic-elliptic type was first considered in \cite{GZ2014}. In \cite{Guo-Zhang}, the authors proved the local exact controllability to a fixed trajectory for Keller-Segel model, where the control acts on the cell density equation.
Moreover, the authors pointed out that the strategy to prove the controllability may not be applied to the case of the control acting on the chemical concentration equation. The main difficulty here is that one cannot obtain the observability estimate for the adjoint equation since
one variable cannot be directly represented by the other in this case.
The local controllability of the Keller-Segel system around a constant trajectory with the control acting on the component of the chemical was discussed in \cite{CG2015}.  The controllability to a constant trajectory referring to the objective trajectory is the constant solution of parabolic-elliptic system, and the controllability to a fixed (non-constant) trajectory denoting the objective trajectory is
the solution of parabolic-parabolic system.
  Later, a controllability result for a chemotaxis-fluid model around some particular trajectories was studied in \cite{CG2017}.
The chemotaxis models in these known results are completely different from our system (\ref{1.2}).
 The strategies in these papers cannot be applied here directly since our system has first order nonlinear couplings. To our knowledge, there is no literature on the controllability to a nonnegative trajectory for systems considered in this paper.

For simplicity, we use notations $L^p(Q)$, $H^p(\O)$ and
$W^{k,p}(\O)$ to denote the $n$ product spaces $L^p(Q)^n$,
$H^p(\O)^n$ and $W^{k,p}(\O)^n$, respectively.
Consider a free system without control function:
\begin{equation}\label{1.5}
\left\{\begin{array}{ll}
\overline{u}_t-\Delta \overline{u}=\nabla\cdot(\overline{u}\overline{\mathbf{v}}),& (x,t)\in Q,\\
\overline{\mathbf{v}}_t-\Delta\overline{\mathbf{v}}=
-\nabla(|\overline{\mathbf{v}}|^2)+\nabla \overline{u},&(x,t)\in Q,\\
\overline{u}=\overline{p}, ~~\overline{\mathbf{v}}=0,&(x,t)\in\Sigma,\\
(\overline{u},\overline{\mathbf{v}})(x,0)=
(\overline{u}_0,\overline{\mathbf{v}}_0)(x),
& x\in\Omega,
\end{array}\right.
\end{equation}
where $(\ou_0,\ov_0)\in H^4(\O)\times H^4(\O)$ satisfies
$$\ou_0-\overline{p}\geq0,\ \mathbf{\ov}_0\geq0\  \text{and}\ \|\ou_0-\overline{p}\|^2_{H^3(\O)}+\|\mathbf{\ov}_0\|^2_{H^3(\O)}\leq \varepsilon$$
for some constant $\varepsilon\in(0,1)$.
Assume that $(\overline{u},\overline{\mathbf{v}})$ is a nonnegative trajectory of equation (\ref{1.5}) associated to  $(\overline{u}_0,\overline{\mathbf{v}}_0)$ and $\overline{p}$. The existence and the nonnegativity of this kind of trajectories will be given in Section 2.
The system (\ref{1.2}) is said to be locally exactly controllable to the trajectory $(\overline{u},\overline{\mathbf{v}})$ at time $T$, if there is a neighborhood $\mathcal{O}$ of $(\overline{u}_0,\overline{\mathbf{v}}_0)$ such that for any initial data $(u_0,\mathbf{v}_0)\in\mathcal{O}$, there exists a control function $\mathbf{h}$  with the corresponding solution $(u,\mathbf{v})$ of (\ref{1.2})  satisfying $$u(x,T)=\overline{u}(x,T),\ \mathbf{v}(x,T)=\overline{\mathbf{v}}(x,T),\ a.e. \ x\in\O.$$

We have the following main result for the system (\ref{1.2}).

\begin{theorem}\label{t4}
Let $r>n+2$, and $(\ou,\ov)$ be the trajectory of system (\ref{1.5}) corresponding to $(\ou_0,\mathbf{\ov}_0)\in H^4(\O)\times H^4(\O)$, which satisfies
\begin{equation*}
\label{2.-1}
\|\ou_0-\overline{p}\|^2_{H^3(\O)}+\|\mathbf{\ov}_0\|^2_{H^3(\O)}\leq \varepsilon,
\ \text{and}\ \ou_0 -\overline{p}\geq 0, ~~\ov_0\geq 0
\end{equation*}
for some constant $\varepsilon\in(0,1)$.
Then, there exists a  constant $\delta>0$, depending only on $n,\o,\O$ and $T$, such that for any $(u_0,\mathbf{v}_0)\in W^{2-\frac{2}{r},r}(\O)\times W^{2-\frac{2}{r},r}(\O)$  satisfying
$$u_0\geq 0,\ \text{and}\ \|u_0-\ou_0\|_{W^{2-\frac{2}{r},r}(\O)}+\|\mathbf{v}_0-\ov_0\|_{W^{2-\frac{2}{r},r}(\O)}\leq \delta,$$
there is a control $\mathbf{h}\in L^r(Q)$, with $\supp\mathbf{h}\subseteq \o\times [0,T]$ and the
 system (\ref{1.2}) satisfies
$$u(x,t)\geq 0\ \text{in}\ Q\ \text{and}\ u(x,T)=\ou(x,T),\
 \mathbf{v}(x,T)=\mathbf{\ov}(x,T)\ \text{in} \ \O.$$
\end{theorem}


\begin{remark}
 In (\ref{1.5}), the nonhomogeneous Dirichlet-Dirichlet boundary conditions ensure that the solution $\ou$ has a positive lower bound. Due to the complexity of system (\ref{1.2}), we need the positive lower bound result of the solution to derive the Carleman estimate.
 Therefore, the boundary condition $\ou=\overline{p}$ is technical, and the strategy developed in this paper cannot be employed to the case of homogeneous Dirichlet-Dirichlet boundary conditions.
\end{remark}

\begin{remark}
We assume that the initial data $(\ou_0,\mathbf{\ov}_0)$ belongs to  $H^4(\O)\times H^4(\O)$ in Theorem \ref{t4}. It is worth mentioning that this regularity of the initial values can be reduced to  $H^3(\O)\times H^3(\O)$ based on the regularizing effect of system (\ref{1.5}), which has been proved in \cite[Theorem 3.1]{T1982}. This strategy has also been used in \cite[Lemma 5]{CMS2020}. It would be quite interesting to study the controllability for more general initial conditions. However, it seems that the method developed in this paper is not enough. We will explain this in Remark \ref{r2}.
\end{remark}

 The study of controllability for coupled parabolic equations has attracted intensive attention in the past few years. In general, the controllability of coupled systems is more difficult than that of single equations.
 Some new phenomena may occur. For instance, the minimal time of control is required for the controllability  of some coupled parabolic systems (see \cite{new}).

 There are many works addressing the controllability of parabolic systems with zero order coupling terms
(see, e.g., \cite{ 2009427, 2009267, 2005426, FAML, 201417, ba, 201091,  2013187} and the rich references therein).
Concerning the case of first order coupling terms,  we refer to  \cite{ yi,  si, ba, shisan, guerrero, shi} and \cite{liu} for some known controllability results for coupled parabolic systems.
In \cite{guerrero}, the author investigated the case of first and second order coupling terms, and the coupling coefficients are constants or only dependent on time variable by means of the Carleman estimate approach.
In \cite{yi} and \cite{shisan}, the controllability
was obtained for some systems with time and space-varying coupling coefficients under some technical conditions. The one-dimensional results were given in \cite{si} and \cite{ba}. Specifically, the main tool in \cite{si} is the moment method and the coefficients only depend on space, while
the authors of \cite{ba} used the fictitious control method to solve the case where the coefficients depend on the space and time variables.
In \cite{shi}, by means of the Lebeau-Robbiano strategy, the internal observability was established for the system with constant or time-dependent coupling terms.
Recently, the work in \cite{liu} studied the case of constant coupling coefficients by an algebraic method.

 Obviously, our system (\ref{1.2}) is a nonlinear coupled parabolic system. The usual way to establish the controllability of a nonlinear system is to prove the controllability of the associated linearized system combined with the fixed point technique.
The key point here is that one needs to establish the suitable observability inequality for the associated adjoint system (see (\ref{1.9})).
To achieve this goal,  we shall employ a similar method as the one used in \cite{guerrero} to  derive a new global Carleman estimate for general coupled parabolic equations.
The main difficulty is the coefficients of the first order coupling terms involving the solutions of the free system (\ref{1.5}), which depend on both time and space variables.
 Hence, it is technically
 more complicated and difficult to deal with  this problem.
 Moreover,  in order to establish the Carleman estimate for the adjoint system, we require that the coupling coefficients belong to $W^{2,1}_\infty(Q)$. Accordingly, another key point of the proof is to show
 the global existence of nonnegative solution for system (\ref{1.5}), and to establish the suitable regularity for the solutions, which has an independent interest even as a pure PDE problem.
 The main idea for obtaining the regularity of the solutions is to use
 the temporal derivatives of the solution to recover the spatial derivatives, since the information of the spatial derivatives of the solution is unknown on the boundary.
 At last, in order to ensure the application of the fixed point argument, we need to improve the regularity of the control, and the technique used to solve this problem is adapted from \cite{barbu2}. Certainly, the existence of first-order coupling terms also makes the proof of this problem more complicated.

The rest of this paper is organized as follows. In Section 2, we investigate the global existence of solutions for system (\ref{1.5}). Section 3 is devoted to showing the null controllability of the linearized system. Then Theorem \ref{t4} is proved in Section 4.

\section{Global existence of the trajectory}

In this section, we first prove the following well-posedness result for system (\ref{1.5}) in order to guarantee
global existence of the trajectory. Then, we will show that the global trajectory preserves the nonnegative property of the initial data.
To the best of our knowledge, the global well-posedness of system (\ref{1.5}) has not been studied in the literature.
It is worth mentioning that the Dirichlet boundary value problem (\ref{1.5})
is meaningful from the biological point of view, see for example \cite{levine2, LZ2015}.

We have the following well-posedness result for (\ref{1.5}).

\begin{theorem}\label{t3}
 If $(\ou_0,\mathbf{\ov}_0)\in H^4(\O)\times H^4(\O)$ satisfies
\begin{equation}
\label{2.-1}
\|\ou_0-\overline{p}\|^2_{H^3(\O)}+\|\mathbf{\ov}_0\|^2_{H^3(\O)}\leq \varepsilon
\end{equation}
for some constant $\varepsilon \in (0,1)$,
 then there exists a unique solution $(\ou,\ov)$ of (\ref{1.5}) satisfying
\bea
&&(\ou,\ov)\in C([0,T];H^4(\O))\times C([0,T];H^4(\O)),\nonumber\\
&&(\ou_t,\ov_t)\in C([0,T];H^1(\O))\cap L^2(0,T; H^2(\O))\times C([0,T];H^1(\O))\cap L^2(0,T; H^2(\O)).\nonumber
\eea
Moreover, it holds that
\begin{equation}
\label{2.-2}
\|\ou-\overline{p}\|^2_{L^\infty(0,T;H^4(\O))}
+\|\ov\|^2_{L^\infty(0,T;H^4(\O))}
\leq C\(\|\ou_0-\overline{p}\|^2_{H^4(\O)}+\|\ov_0\|^2_{H^4(\O)}\).
\end{equation}
\end{theorem}

Let $\overline{w}=\overline{u}-\overline{p}$. We will start by studying the following system of $(\overline{w},\ov)$:
\begin{equation}\label{2.00}
\left\{\begin{array}{ll}
\overline{w}_t-\Delta \overline{w}=\nabla\cdot(\overline{w}\overline{\mathbf{v}})+\overline{p}\nabla\cdot\overline{\mathbf{v}},& (x,t)\in Q,\\
\overline{\mathbf{v}}_t-\Delta\overline{\mathbf{v}}=
-\nabla(|\overline{\mathbf{v}}|^2)+\nabla \overline{w},&(x,t)\in Q,\\
\overline{w}=\overline{\mathbf{v}}=0,&(x,t)\in\Sigma,\\
(\overline{w},\overline{\mathbf{v}})(x,0)=
(\overline{w}_0,\overline{\mathbf{v}}_0)(x)=(\overline{u}_0-\overline{p},\overline{\mathbf{v}}_0)(x),
& x\in\Omega.
\end{array}\right.
\end{equation}

The proof of Theorem \ref{t3} is based on the standard continuity argument. Hence,  we first need  to assume that
there exists a small positive constant $\delta_0<1$ satisfying
\begin{equation}
\label{2.0}
\sup_{0\leq t\leq T}(\|\ow\|_{H^3(\Omega)}^2+\|\mathbf{\ov}\|_{H^3(\Omega)}^2) <\delta_0.
\end{equation}
In what follows, we will establish some \emph{a priori} estimates to close (\ref{2.0}).

\begin{lemma}\label{l20}
Under the assumption (\ref{2.0}), it holds that
\begin{equation}
\label{2.1}
\sup_{0\leq t\leq T}( \|\ow\|_{H^1(\Omega)}^2+ \|\mathbf{\ov}\|_{H^1(\Omega)}^2)+\int^T_0 \big(\|\nabla\ow\|_{H^1(\Omega)}^2+ \|\nabla\mathbf{\ov}\|_{H^1(\Omega)}^2\big)dt
 \leq C(\|\ow_0\|_{H^1(\Omega)}^2+ \|\mathbf{\ov}_0\|_{H^1(\Omega)}^2).
\end{equation}
\end{lemma}

\no {\bf Proof.}
Multiplying (\ref{2.00}) by $\ow$ and $\overline{p}\mathbf{\ov}$ respectively,  integrating over $\O$,  and using integration by parts,
H\"{o}lder, Poincar\'{e}, Young and Cauchy inequalities, we have
\bea\label{2.2}
&&\frac{1}{2}\frac{d}{dt}\int_{\O}(\ow^2+\overline{p}\ov^2)dx+
\int_\O(|\nabla\ow|^2+ \overline{p}|\nabla\ov|^2)dx
\nonumber\\[2mm]
&&
\leq C\|\nabla\ow\|_{L^2(\O)}\|\mathbf{\ov}\|_{L^2(\O)}\|\ow\|_{L^\infty(\O)}+C\|\nabla\ov\|_{L^2(\O)}^2
\|\ov\|_{L^\infty(\O)}
\nonumber\\[2mm]
&&
\leq C\|\nabla\ow\|_{L^2(\O)}\|\mathbf{\ov}\|_{L^2(\O)}\|\ow\|_{H^2(\O)}+C\|\nabla\ov\|_{L^2(\O)}^2
\|\ov\|_{H^2(\O)}
\nonumber\\[2mm]
&&
\leq C\delta_0(\|\nabla\ow\|^2_{L^2(\O)}+ \|\nabla\mathbf{\ov}\|^2_{L^2(\O)}),
\eea
where the assumption (\ref{2.0}) is used in the last inequality.

Then, taking the $L^2$ inner product of the first equation in $(\ref{2.00})$ with $\Delta\ow$ and the second equation in $(\ref{2.00})$ with $\overline{p}\Delta\ov$, we derive
\bea
\label{2.5}
&&\frac{1}{2}\frac{d}{dt}\int_{\O}|\nabla\ow|^2dx+
\int_\O|\nabla^2\ow|^2dx\nonumber\\[2mm]
&&\leq C\|\nabla^2\ow\|_{L^2(\O)}(\|\nabla\ow\|_{L^2(\O)}\|\ov\|_{L^\infty(\O)} +\|\nabla\ov\|_{L^2(\O)}\|\ow\|_{L^\infty(\O)}) + \overline{p}\|\na^2\ow\|_{L^2(\O)}\|\na\ov\|_{L^2(\O)}
\nonumber\\[2mm]
&&\leq C\|\nabla^2\ow\|_{L^2(\O)}(\|\nabla\ow\|_{L^2(\O)}\|\ov\|_{H^2(\O)} + \|\nabla\ov\|_{L^2(\O)}\|\ow\|_{H^2(\O)})+\overline{p}\|\na^2\ow\|_{L^2(\O)}\|\na\ov\|_{L^2(\O)}\nonumber\\[2mm]
&&\leq C \delta_0(\|\nabla^2\ow\|^2_{L^2(\O)}+\|\nabla\ow\|^2_{L^2(\O)}+ \|\nabla\mathbf{\ov}\|^2_{L^2(\O)})
+\frac{\overline{p}^2}{3}\|\na\ov\|^2_{L^2(\O)}+ \frac{3}{4}\|\na^2\ow\|^2_{L^2(\O)}
\eea
and
\bea
\label{2.6}
&&\frac{1}{2}\frac{d}{dt}\int_\O\overline{p}|\na\ov|^2dx+
\int_\O\overline{p}|\na^2\ov|^2dx\nonumber\\
&&\leq C\|\nabla^2\ov\|_{L^2(\O)}\|\nabla\ov\|_{L^2(\O)}\|\ov\|_{L^\infty(\O)}+ \overline{p}\|\na\ow\|_{L^2(\O)}\|\na^2\ov\|_{L^2(\O)}
\nonumber\\
&& \leq C \delta_0(\|\nabla^2\ov\|^2_{L^2(\O)}+ \|\nabla\ov\|^2_{L^2(\O)})+
\frac{\overline{p}}{3}\|\na\ow\|^2_{L^2(\O)}+\frac{3\overline{p}}{4}
\|\na^2\ov\|^2_{L^2(\O)}.
\eea
It follows from (\ref{2.5}) and (\ref{2.6}) with $\delta_0$ small enough that
\bea\label{2.7}
&&\frac{1}{2}\frac{d}{dt}\int_\O\(|\na\ow|^2+\overline{p}|\na\ov|^2\)dx
+C_1\int_\O(|\na^2\ow|^2+\overline{p}|\na^2\ov|^2)dx\nonumber\\
&&\leq C\delta_0(\|\nabla\ow\|^2_{L^2(\O)}+ \|\nabla\mathbf{\ov}\|^2_{L^2(\O)})
+\frac{\overline{p}}{3}\|\na\ow\|^2_{L^2(\O)}+\frac{\overline{p}^2}{3}\|\na\ov\|^2_{L^2(\O)},
\eea
where $C_1>0$.

Multiplying (\ref{2.2}) by $\overline{p}$, adding the resulting inequality with (\ref{2.7}), and integrating over $[0,t]$, we obtain (\ref{2.1}) for $\delta_0$ small enough.
The proof is completed.
\endpf

Obviously, we can get the following corollary by applying
Lemma \ref{l20}  directly to system (\ref{2.00}).

\begin{corollary}\label{c21}
Under the assumption (\ref{2.0}), it holds that
\begin{equation}
\label{2.8}
\int^T_0 \big(\|\ow_t\|_{L^2(\Omega)}^2+ \|\ov_t\|_{L^2(\Omega)}^2\big)dt
 \leq C(\|\ow_0\|_{H^1(\Omega)}^2+ \|\mathbf{\ov}_0\|_{H^1(\Omega)}^2).
\end{equation}
\end{corollary}

Next we shall turn to the estimation of higher order spatial derivatives of the solution.
Because of the lack of information of the spatial derivatives of the solution on the boundary,
we need to use temporal derivatives and system (\ref{2.00}) to obtain bounds for the spatial derivatives.

\begin{lemma}\label{l22}
Under the assumption (\ref{2.0}), it holds that
\begin{eqnarray}
&&\label{2.9}
\sup_{0\leq t\leq T}(\|\ow_t\|_{H^1(\Omega)}^2+ \|\ov_t\|_{H^1(\Omega)}^2)+\int^T_0 \big(\|\nabla\ow_t\|_{H^1(\Omega)}^2+ \|\nabla\mathbf{\ov}_t\|_{H^1(\Omega)}^2\big)dt\nonumber\\[2mm]
&& \leq C(\|\ow_0\|_{H^3(\Omega)}^2+ \|\mathbf{\ov}_0\|_{H^3(\Omega)}^2).
\end{eqnarray}
\end{lemma}

\no {\bf Proof.}
Differentiating the first equation of $(\ref{2.00})$ with respect to $t$, multiplying the resulting equation
by $\ow_t$ and integrating over $\Omega$, we have
\bea
&&\frac{1}{2}\frac{d}{dt}\int_{\O}\ow_t^2dx+
\int_\O|\nabla\ow_t|^2dx=-\int_{\O}\nabla\ow_t \cdot(\ow\ov)_t dx - \overline{p}\int_{\O}\nabla\ow_t \cdot\ov_t dx\nonumber\\
&&\leq \frac{1}{2}\|\nabla\ow_t\|_{L^2(\O)}^2 + C(\|\ow_t\|_{L^2(\O)}^2\|\ov\|_{L^\infty(\O)}^2 + \|\ov_t\|_{L^2(\O)}^2\|\ow\|_{L^\infty(\O)}^2+\|\ov_t\|_{L^2(\O)}^2),\nonumber
\eea
which implies
\bea \label{2.10}
\frac{d}{dt}\int_{\O}\ow_t^2dx+
\int_\O|\nabla\ow_t|^2dx
\leq C\delta_0(\|\ow_t\|_{L^2(\O)}^2  + \|\ov_t\|_{L^2(\O)}^2)+ C\|\ov_t\|_{L^2(\O)}^2.
\eea
Moreover,
it follows from equations $(\ref{2.00})$ that for any $0\leq t\leq T$,
\bea \label{2.11}
\|\ow_t(t)\|_{L^2(\O)}^2
\leq C (1+ \|\ov(t)\|_{H^2(\O)}^2) (\|\ow(t)\|_{H^2(\O)}^2+\|\ov(t)\|_{H^1(\O)}^2)
\eea
and
\bea \label{2.12}
\|\ov_t(t)\|_{L^2(\O)}^2
\leq C (1+ \|\ov(t)\|_{H^2(\O)}^2) (\|\ov(t)\|_{H^2(\O)}^2+\|\ow(t)\|_{H^1(\O)}^2).
\eea
Thus, integrating (\ref{2.10}) over $[0, t]$ and using (\ref{2.8}),
we obtain
\bea \label{2.13}
&&\sup_{0\leq t\leq T}\|\ow_t\|_{L^2(\Omega)}^2 +\int^T_0  \|\nabla\ow_t\|_{L^2(\Omega)}^2 dt\nonumber\\
&&\leq \|\ow_t(0)\|_{L^2(\Omega)}^2+ C\int^T_0(\|\ow_t\|_{L^2(\O)}^2  + \|\ov_t\|_{L^2(\O)}^2)dt\\
&&\leq C(\|\ow_0\|_{H^2(\Omega)}^2+ \|\mathbf{\ov}_0\|_{H^2(\Omega)}^2),\nonumber
 \eea
where we have used the fact from (\ref{2.11}) that
$\|\ow_t(0)\|_{L^2(\O)}^2
\leq C (\|\ow_0\|_{H^2(\Omega)}^2+ \|\mathbf{\ov}_0\|_{H^2(\Omega)}^2)$.
Similarly, with the help of (\ref{2.8}) and (\ref{2.12}), for $\ov$, it holds that
\bea \label{2.14}
\sup_{0\leq t\leq T}\|\ov_t\|_{L^2(\Omega)}^2 +\int^T_0  \|\nabla\ov_t\|_{L^2(\Omega)}^2 dt
 \leq C(\|\ow_0\|_{H^2(\Omega)}^2+ \|\mathbf{\ov}_0\|_{H^2(\Omega)}^2).
 \eea

Next, we take $\frac{\partial}{\partial t}$ to (\ref{2.00}), multiply the resulting equations by $\Delta\ow_t$ and $\Delta\ov_t$ respectively and use  integration by
parts to derive
\bea
&&\frac{1}{2}\frac{d}{dt}\int_{\O}\big(|\nabla\ow_t|^2+ |\nabla\ov_t|^2\big)dx+
\int_\O \big(|\nabla^2\ow_t|^2+ |\nabla^2\ov_t|^2\big)dx\nonumber\\
&&= -\int_\O (\nabla(\ow\ov))_t\cdot \Delta \ow_t dx - \overline{p}\int_\O \nabla\cdot\ov_t\cdot \Delta \ow_t dx
+\int_\O (\nabla(|\ov|^2))_t\cdot \Delta \ov_t dx -\int_\O \nabla\ow_t\cdot \Delta \ov_t dx
\nonumber\\
&&\leq \frac{1}{2}\!\int_\O \!\big(|\nabla^2\ow_t|^2+ |\nabla^2\ov_t|^2\big)dx
\!+\!C\left(\|\ov\|_{L^\infty(\Omega)}^2\int_\O \big(|\nabla\ow_t|^2+ |\nabla\ov_t|^2\big)dx
\!+\!\|\ow\|_{L^\infty(\Omega)}^2 \int_\O |\nabla\ov_t|^2 dx\right.
\nonumber\\
&&
\left.\quad+
\|\nabla\ov\|_{L^\infty(\Omega)}^2\int_\O \big(|\ow_t|^2+ |\ov_t|^2\big)dx
+ \|\nabla\ow\|_{L^\infty(\Omega)}^2  \int_\O |\ov_t|^2dx  +\int_\O |\nabla\ov_t|^2 dx+\int_\O |\nabla\ow_t|^2 dx\right),\nonumber
\eea
which implies
\bea \label{2.15}
&&\frac{d}{dt}\int_{\O}\big(|\nabla\ow_t|^2+ |\nabla\ov_t|^2\big)dx+
\int_\O \big(|\nabla^2\ow_t|^2+ |\nabla^2\ov_t|^2\big)dx\nonumber\\[2mm]
&&\leq C\int_\O \big(|\nabla\ow_t|^2+ |\nabla\ov_t|^2+ |\ow_t|^2+ |\ov_t|^2\big)dx.
 \eea
Taking the derivative of $(\ref{2.00})$ with respect to $x$, by some straightforward calculations, we have, for any $0\leq t\leq T$,
\bea \label{2.16}
\|\nabla\ow_t(t)\|_{L^2(\O)}^2
\leq C (1+ \|\ov(t)\|_{H^3(\O)}^2)( \|\ow(t)\|_{H^3(\O)}^2+\|\ov(t)\|_{H^2(\O)}^2)
\eea
and
\bea \label{2.17}
\|\nabla\ov_t(t)\|_{L^2(\O)}^2
\leq C (1+ \|\ov(t)\|_{H^3(\O)}^2) (\|\ov(t)\|_{H^3(\O)}^2+\|\ow(t)\|_{H^2(\O)}^2).
\eea
Now, integrating (\ref{2.15}) over $[0,t]$, by (\ref{2.8}), (\ref{2.13}), (\ref{2.14}), (\ref{2.16}) and (\ref{2.17}), we arrive at
\bea \label{2.18}
&&\sup_{0\leq t\leq T}(\|\nabla\ow_t\|_{L^2(\Omega)}^2+ \|\nabla\ov_t\|_{L^2(\Omega)}^2) +\int^T_0 \big(\|\nabla^2\ow_t\|_{L^2(\Omega)}^2+ \|\nabla^2\ov_t\|_{L^2(\Omega)}^2 \big)dt\nonumber\\[2mm]
 &&\leq C(\|\ow_0\|_{H^3(\Omega)}^2+ \|\mathbf{\ov}_0\|_{H^3(\Omega)}^2).
\eea
This together with (\ref{2.13}) and (\ref{2.14}) yields (\ref{2.9}).
The proof is completed.
\endpf

With these \emph{a priori} estimates at hand, we are ready to close the assumption (\ref{2.0}).

\no{\bf Proof of Theorem \ref{t3}.}
It follows from the equations (\ref{2.00}) and all the estimates above that
\bea \label{2.19}\sup_{0\leq t\leq T}(\|\ow\|^2_{H^3(\O)}+\|\ov\|^2_{H^3(\O)})
\leq C(\|\ow_0\|^2_{H^3(\O)}+\|\ov_0\|^2_{H^3(\O)})\leq C\varepsilon.\eea
If $\varepsilon$ is suitably small such that $C\varepsilon\leq \delta_0$, by the standard continuity argument (see \cite{MN1980,TY2018}), the
estimate (\ref{2.0}) is closed.
Notice that the local existence and uniqueness of the solution to the equations (2.3) can be established by using the classical theory of linear parabolic system (see, for example, \cite[p.616]{LSU1968}) combining with Schauder fixed point theorem.
Thus, applying (\ref{2.0}) and all \emph{a priori} estimates, we extend the local solutions to be a global
solution and the uniqueness of global solution in $C([0,T];H^3(\O))$ is guaranteed by the uniqueness
of local solution.

Thus, it only remains to establish the regularity in $H^4$ space.
In view of equations $(\ref{2.00})$, $(\ref{2.16})$ and $(\ref{2.17})$, it holds that for any $0\leq t\leq T$,
\bea
&\|\ow_{tt}(t)\|_{L^2(\O)}^2
\leq C (1+ \|\ov(t)\|_{H^4(\O)}^2) (\|\ow(t)\|_{H^4(\O)}^2+\|\ov(t)\|_{H^3(\O)}^2), \label{2.21}\\
&\|\ov_{tt}(t)\|_{L^2(\O)}^2
\leq C (1+ \|\ov(t)\|_{H^4(\O)}^2) (\|\ov(t)\|_{H^4(\O)}^2+\|\ow(t)\|_{H^3(\O)}^2), \label{2.22}
\eea
and
\bea
&\|\ow_{tt}(t)\|_{L^2(\O)}^2\leq C \(\|\nabla^2\ow_t \|_{L^2(\O)}^2+\|\ow_t \|_{H^1(\O)}^2\|\ov \|_{H^3(\O)}^2+\|\ov_t \|_{H^1(\O)}^2(1+\|\ow \|_{H^3(\O)}^2 )\),\label{2.23}\\
&\|\ov_{tt}(t)\|_{L^2(\O)}^2\leq C \big(\|\nabla^2\ov_t \|_{L^2(\O)}^2+\|\ov_t \|_{H^1(\O)}^2\|\ov \|_{H^3(\O)}^2+ \|\nabla\ow_t \|_{L^2(\O)}^2\big). \label{2.24}
\eea

Differentiating the first equation of $(\ref{2.00})$ with respect to $t$ twice, and then
multiplying the resulting equation
by $\ow_{tt}$ and integrating over $\Omega$, we have
\bea
&&\frac{1}{2}\frac{d}{dt}\int_{\O}\ow_{tt}^2dx+
\int_\O|\nabla\ow_{tt}|^2dx\nonumber=-\int_{\O}\nabla\ow_{tt} \cdot(\ow\ov)_{tt} dx - \overline{p}\int_{\O}\nabla\ow_{tt} \cdot\ov_{tt} dx\nonumber\\[2mm]
&&\leq \frac{1}{2}\|\nabla\ow_{tt}\|_{L^2(\O)}^2 + C(\|\ow_{tt}\|_{L^2(\O)}^2\|\ov\|_{L^\infty(\O)}^2+ \|\ov_{tt}\|_{L^2(\O)}^2\|\ow\|_{L^\infty(\O)}^2+ \|\ov_{t}\|_{L^2(\O)}^2\|\ow_t\|_{L^\infty(\O)}^2\nonumber\\[2mm]
&&\quad+\|\ov_{tt}\|_{L^2(\O)}^2),
\eea
which together with $(\ref{2.14})$ and $(\ref{2.19})$ implies
\bea\label{2.25}
&&\frac{d}{dt}\int_{\O}\ow_{tt}^2dx+
\int_\O|\nabla\ow_{tt}|^2dx\nonumber\\
&&\leq C(\|\ow_{tt}\|_{L^2(\O)}^2\|\ov\|_{H^2(\O)}^2 + \|\ov_{tt}\|_{L^2(\O)}^2\|\ow\|_{H^2(\O)}^2 + \|\ov_{t}\|_{L^2(\O)}^2\|\ow_t\|_{H^2(\O)}^2+\|\ov_{tt}\|_{L^2(\O)}^2)
\nonumber\\
&&\leq
C(\|\ow_0\|^2_{H^3(\O)}+\|\ov_0\|^2_{H^3(\O)})(\|\ow_{tt}\|_{L^2(\O)}^2+ \|\ov_{tt}\|_{L^2(\O)}^2 + \|\ow_t\|_{H^2(\O)}^2)+C\|\ov_{tt}\|_{L^2(\O)}^2.
\eea
Then, integrating (\ref{2.25}) over $[0,t]$, by
$(\ref{2.8})$, $(\ref{2.13})$, $(\ref{2.14})$, $(\ref{2.18})$ and $(\ref{2.21})$-$(\ref{2.24})$, we obtain
\bea \label{2.26}
\sup_{0\leq t\leq T}\|\ow_{tt}\|_{L^2(\Omega)}^2+\int^T_0 \|\nabla\ow_{tt}\|_{L^2(\Omega)}^2dt
 \leq C(\|\ow_0\|_{H^4(\Omega)}^2+ \|\mathbf{\ov}_0\|_{H^4(\Omega)}^2).
\eea
Similarly, for $\ov$, it holds that
\bea \label{2.27}
\sup_{0\leq t\leq T}\|\ov_{tt}\|_{L^2(\Omega)}^2+\int^T_0 \|\nabla\ov_{tt}\|_{L^2(\Omega)}^2dt
 \leq C(\|\ow_0\|_{H^4(\Omega)}^2+ \|\mathbf{\ov}_0\|_{H^4(\Omega)}^2).
\eea
Therefore, with the help of the equations $(\ref{2.00})$, $(\ref{2.19})$, $(\ref{2.26})$ and $(\ref{2.27})$,
we deduce $(\ref{2.-2})$.
This completes the proof of Theorem \ref{t3}.
\endpf

\begin{corollary}\label{c22}
 Assume that the conditions in Theorem \ref{t3} hold, and   $\ou_0-\overline{p}\geq 0$, $\ov_0\geq 0$.
Then the solution of (\ref{1.5}) satisfies
$$\ou -\overline{p}\geq 0, ~~\ov\geq 0,\ \forall\ (x,t)\in Q.$$
\end{corollary}
\no {\bf Proof.}
It follows from  Theorem \ref{t3} that the global trajectory $(\ou, \ov)$ is the classical solution of
(\ref{1.5}). Obviously, $(\overline{p}, 0)$ can be regarded as a lower solution of (\ref{1.5}). Thus, the conclusion
of this corollary follows from the comparison principle immediately.
\endpf

\section{Null controllability of the linearized system}

Let $y=u-\overline{u}$, $\mathbf{z}=\mathbf{v}-\overline{\mathbf{v}}$, $y_0=u_0-\overline{u}_0$ and $\mathbf{z}_0=\mathbf{v}_0-\overline{\mathbf{v}}_0$.
An easy computation shows that
$(y,\mathbf{z})$ satisfies
\begin{equation}\label{1.6}
\left\{\begin{array}{ll}
y_t-\Delta y=\nabla\cdot(y(\mathbf{z}+\overline{\mathbf{v}}))
+\nabla\cdot(\overline{u}\mathbf{z}),& (x,t)\in Q,\\
\mathbf{z}_t-\Delta\mathbf{z}=
-\nabla(|\mathbf{z}|^2+2\overline{\mathbf{v}}\cdot\mathbf{z})
+\nabla y+
\chi_{\o}\mathbf{h}, &(x,t)\in Q,\\ y=\mathbf{z}=0,&(x,t)\in\Sigma,\\
(y,\mathbf{z})(x,0)=(y_0,\mathbf{z}_0)(x),
& x\in\Omega.
\end{array}\right.
\end{equation}
Obviously, the local exact controllability to the trajectory $(\overline{u},\overline{\mathbf{v}})$ for  equations (\ref{1.2}) is equivalent to the local  null controllability  of system (\ref{1.6}).

 For $p\geq 2,$ define the Banach space  $V^p$  by
 $$V^p:=\{y: y\in L^p(0,T;W^{2,p}(\O)\cap W^{1,p}_0(\O));y_t\in L^p(Q)\},$$
 and  its natural norm  $\|\cdot\|_{V^p}$ by $\|y\|_{V^p}
 =\|y\|_{L^p(0,T;W^{2,p}(\O))}+\|y_t\|_{L^p(Q)}$.

In this section,  we consider the null controllability of the following linearized system of (\ref{1.6}):
\begin{equation}\label{1.7}
\left\{\begin{array}{ll}
y_t-\Delta y=\nabla\cdot(\mathbf{a}y)
+\nabla\cdot(B\mathbf{z}),& (x,t)\in Q,\\
\mathbf{z}_t-\Delta\mathbf{z}=
-\nabla(\mathbf{b}\cdot\mathbf{z})
+\nabla y+
\chi_{\o}\mathbf{h}, &(x,t)\in Q,\\
y=\mathbf{z}=0,&(x,t)\in\Sigma,\\
(y,\mathbf{z})(x,0)=(y_0,\mathbf{z}_0)(x),
& x\in\Omega,
\end{array}\right.
\end{equation}
where $\mathbf{h}$ is the control, $(y_0,\mathbf{z}_0)$
 is the given initial value, and
\begin{equation}\label{e82}
\mathbf{a}, \mathbf{b},
\nabla\mathbf{a}, \nabla\mathbf{b}\in  L^\infty(Q),\
 B, B_t, \na B, \Delta B\in L^\infty(Q)\  \text{and}\  B\  \text{has  a positive lower bound}.
 \end{equation}
 In fact, $\mathbf{a}=\bm{\eta}+\overline{\mathbf{v}}$,
$\mathbf{b}=\bm{\eta}+2\overline{\mathbf{v}}$ and $B=\overline{u}$, where $\bm{\eta}\in  V^r$ is a known function, and $r>n+2$.

Write
$$M_1=1+T\(1+\|\na\cdot\mathbf{a}\|_{L^\infty(Q)}
\!+\|\na\cdot\mathbf{b}\|_{L^\infty(Q)}
\!+\|\mathbf{a}\|^2_{L^\infty(Q)}
\!+\|\mathbf{b}\|^2_{L^\infty(Q)}
\!+\|B\|^2_{L^\infty(Q)}+\|\na B\|_{L^\infty(Q)}\),$$
$$M_2=1+\|\na\cdot\mathbf{a}\|_{L^\infty(Q)}
+\|\na\cdot\mathbf{b}\|_{L^\infty(Q)}
+\|\mathbf{a}\|_{L^\infty(Q)}
+\|\mathbf{b}\|_{L^\infty(Q)}
+\|B\|_{L^\infty(Q)}+\|\na B\|_{L^\infty(Q)}.$$
 We have the following  well-posedness result for system (\ref{1.7}).

\begin{proposition}\label{p1}
Assume that $\mathbf{a}, \mathbf{b}, B\in L^\infty(Q),\
\nabla\cdot\mathbf{a}, \nabla\cdot\mathbf{b}, \na B\in  L^\infty(Q),$ $y_0,\mathbf{z}_0\in W^{2-\frac{2}{p},p}(\O)\cap H_0^1(\O)$,
  and $\mathbf{h}\in L^p(\o\times (0,T))$ with $p\geq2$ being arbitrary. Then system (\ref{1.7}) admits a unique strong solution
$(y,\mathbf{z})\in V^p\times V^p$. Moreover,
there exist positive constants $C=C(\O,n,p)$ and $k_1=k_1(n)$ such that
\bel{4.1}
\|(y,\mathbf{z})\|_{V^p\times V^p}
\leq e^{CM_1}M_2^{k_1}\(\|(y_0,\mathbf{z}_0)\|
_{W^{2-\frac{2}{p},p}(\O)\times W^{2-\frac{2}{p},p}(\O)}
+\|\mathbf{h}\|_{L^p(\o\times (0,T))}\).
\ee
\end{proposition}

\no {\bf Proof.}
 We split the proof into two steps.

 \no {\bf Step 1.} When $p=2$, multiplying  the first equation of (\ref{1.7}) by $y$ and integrating it on $\O$, we get
 \begin{eqnarray*}
 &&\frac{1}{2}\frac{d}{dt}\int_\O y^2dx+
 \int_\O|\na y|^2dx
 =\int_\O \mathbf{a}y\cdot\na ydx+\int_\O y^2\na\cdot \mathbf{a}dx
 -\int_\O B \mathbf{z}\cdot\na y dx\\
 &&\leq \(\|\mathbf{a}\|_{L^{\infty}(\O)}^2+\|\na\cdot \mathbf{a}\|_{L^{\infty}(\O)}\)\int_\O |y|^2dx
 +\frac{1}{2}\int_\O|\na y|^2dx
 +\|B\|_{L^\infty(\O)}^2\int_\O|\mathbf{z}|^2dx.
 \end{eqnarray*}
Doing the same thing to the second equation of (\ref{1.7}), we obtain
 \begin{eqnarray*}
&&\frac{1}{2}\frac{d}{dt}\int_\O \mathbf{z}^2dx+
 \int_\O|\na \mathbf{z}|^2dx
 =\int_\O \mathbf{b}\cdot\mathbf{z}\na \cdot\mathbf{z}dx
 +\int_\O \mathbf{z}\cdot\na y dx+\int_\O \chi_{\o}\mathbf{h}\cdot\mathbf{z} dx\\
 &&\leq C\(1+\|\mathbf{b}\|_{L^{\infty}(\O)}^2\)\int_\O |\mathbf{z}|^2dx
 +\frac{1}{4}\int_\O|\na \mathbf{z}|^2dx
 +\frac{1}{4}\int_\O|\na y|^2dx
 +\int_\O|\chi_{\o}\mathbf{h}|^2dx.
 \end{eqnarray*}
Then,
\begin{eqnarray*}
 &&\frac{1}{2}\frac{d}{dt}\int_\O \(y^2+\mathbf{z}^2\)dx+
 \int_\O\(|\na y|^2+|\na \mathbf{z}|^2\)dx\\
&&\leq C\(\|\mathbf{a}\|_{L^{\infty}(\O)}^2
+\|\nabla\cdot\mathbf{a}\|_{L^{\infty}(\O)}\)\int_\O |y|^2dx
 +C\(1+\|\mathbf{b}\|_{L^{\infty}(\O)}^2+\|B\|_{L^\infty(\O)}^2\)\int_\O |\mathbf{z}|^2dx\\[2mm]
&&\quad+C\int_\o|\mathbf{h}|^2dx.
 \end{eqnarray*}
By Gronwall's inequality, we have
\bea\label{4.2}
&&\int_\O\(y^2(x,t)+\mathbf{z}^2(x,t)\)dx
+\int_0^t\int_\O\(|\na y|^2+|\na \mathbf{z}|^2\)dxdt\nonumber\\[2mm]
&&\leq Ce^{CM_1}\(\int_\O  (y_0^2+\mathbf{z}_0^2)dx
+\int_0^T\int_\o |\mathbf{h}|^2dxdt\).
\eea

On the other hand, by the first equation of (\ref{1.7}), we have
$$(y_t-\Delta y)^2=\(\na\cdot(\mathbf{a}y)+\na\cdot(B\mathbf{z})\)^2.$$
Integrating the previous inequality  on $\O\times(0,t)$, we obtain
\begin{eqnarray*}
&&\int_0^t\int_\O y_t^2dxdt
+\int_0^t\int_\O |\Delta y|^2 dxdt
-2\int_0^t\int_\O y_t\Delta ydxdt\nonumber\\
&&=\int_0^t\int_\O y_t^2dxdt
+\int_0^t\int_\O |\Delta y|^2 dxdt
+\int_0^t\int_\O\frac{d}{dt}|\na y|^2dxdt
\nonumber\\
&&=\int_0^t\int_\O \(\na\cdot(\mathbf{a}y)+\na\cdot(B\mathbf{z})\)^2dxdt.
\end{eqnarray*}
Combining this with  (\ref{4.2}), we have
\begin{eqnarray}\label{4.6}
&&\int_0^t\int_\O y_t^2dxdt
+\int_0^t\int_\O |\Delta y|^2 dxdt
+\int_\O |\na y(x,t)|^2dx\nonumber\\
&&=\int_0^t\int_\O \(\na\cdot(\mathbf{a}y)+\na\cdot(B\mathbf{z})\)^2dxdt+\int_\O|\na y_0|^2dx\nonumber\\
&&\leq Ce^{CM_1}M_2^{k_1}\[\int_\O\(y_0^2+\mathbf{z}_0^2+|\na y_0|^2\)dx
+\int_0^T\int_\o |\mathbf{h}|^2dxdt\],
\end{eqnarray}
where $k_1=k_1(n),$ $C=C(n,\O,p)$.
 Similarly, we deal with the second equation of (\ref{1.7}), which implies
\begin{eqnarray}\label{4.5}
&&\int_0^t\int_\O \mathbf{z}_t^2dxdt
+\int_0^t\int_\O |\Delta \mathbf{z}|^2 dxdt
+\int_\O|\na \mathbf{z}(x,t)|^2 dx\nonumber\\
&& \leq Ce^{CM_1}M_2^{k_1}\[\int_\O\(y_0^2+\mathbf{z}_0^2+|\na \mathbf{z}_0|^2\)dx
+\int_0^T\int_\o |\mathbf{h}|^2dxdt\].
\end{eqnarray}
By (\ref{4.6}) and (\ref{4.5}), we have
\begin{eqnarray*}
&&\int_0^t\int_\O \(|y_t|^2+|\mathbf{z}_t|^2\)dxdt
+\int_0^t\int_\O \(|\Delta y|^2 +|\Delta \mathbf{z}|^2\) dxdt
+\int_\O\(|\na y(x,t)|^2+|\na \mathbf{z}(x,t)|^2\) dx\nonumber\\
&&\leq Ce^{CM_1}M_2^{k_1}\[\int_\O\(y_0^2+\mathbf{z}_0^2+|\na y_0|^2+|\na \mathbf{z}_0|^2\)dx
+\int_0^T\int_\o |\mathbf{h}|^2dxdt\].
\end{eqnarray*}

 \no {\bf Step 2.} We consider the case $p>2$.
 We only show the case when $n=3$, since the proof is similar when $n=1$ or $2$.

 By Step 1, we know that the solution of (\ref{1.7}) lies in $V^2\times V^2$, and
\bel{e4.1}
\|(y,\mathbf{z})\|_{V^2\times V^2}
\leq e^{CM_1}M_2^{k_1}\(\|(y_0,\mathbf{z}_0)\|
_{W^{1,2}(\O)\times W^{1,2}(\O)}
+\|\mathbf{h}\|_{L^2(\o\times (0,T))}\).
\ee

 On the one hand, let $\mathbf{f}_1=-\mathbf{z}\nabla\cdot\mathbf{b}
-\mathbf{b}\nabla\cdot\mathbf{z}
+\nabla y+\chi_\o\mathbf{h}$.
 Note that $\mathbf{b},\na\cdot\mathbf{b}\in L^\infty(Q)$, by Sobolev embedding theorem,  we have
 $$\na y, \mathbf{b}\na\cdot\mathbf{z}\in L^2(0,T;L^q(\O))\cap
 L^\infty(0,T;L^2(\O)),$$
 where $q=\frac{2n}{n-2}$.
 Therefore, $\mathbf{f_1}\in L^p(0,T;L^{p_1}(\O))$, with
$p_1=\min\{p,\frac{2np}{np-4}\},$
and
\bel{e4.2}
\|\mathbf{f_1}\|_{L^p(0,T;L^{p_1}(\O))}
\leq C\(1+\|\mathbf{b}\|_{L^\infty(Q)}
+\|\na\cdot\mathbf{b}\|_{L^\infty(Q)}\)\cdot
\(\|(y,\mathbf{z})\|_{V^2\times V^2}+\|\mathbf{h}\|_{L^p(\o\times (0,T))}\).
\ee
 Then by Theorem 2.3 in \cite{2006-18}, we obtain
 $$\mathbf{z}\in L^p(0,T;W^{2,p_1}(\O)),\ \mathbf{z}_t\in L^p(0,T;L^{p_1}(\O)),$$
 and
 \begin{equation*}
 \|\mathbf{z}\|_{L^p(0,T;W^{2,p_1}(\O))}
 +\|\mathbf{z}_t\|_{L^p(0,T;L^{p_1}(\O))}
 \leq C\(\|\mathbf{f_1}\|_{L^p(0,T;L^{p_1}(\O))}
 +\|\mathbf{z}_0\|_{W^{2-\frac{2}{p},p}(\O)}\),
 \end{equation*}
 where $C>0$ is a constant independent of $T$.
 Combining this with (\ref{e4.2}), we get
 \begin{eqnarray}\label{e4.3}
 &&\|\mathbf{z}\|_{L^p(0,T;W^{2,p_1}(\O))}
 +\|\mathbf{z}_t\|_{L^p(0,T;L^{p_1}(\O))}\nonumber\\
 &&\leq \!C\(1\!+\!\|\mathbf{b}\|_{L^\infty(Q)}
\!+\!\|\na\!\cdot\!\mathbf{b}\|_{L^\infty(Q)}\)
\(\|(y,\mathbf{z})\|_{V^2\times V^2}\!+\!\|\mathbf{h}\|_{L^p(\o\times (0,T))}
\!+\!\|\mathbf{z}_0\|_{W^{2-\frac{2}{p},p}(\O)}\).
 \end{eqnarray}

 On the other hand,
 take $f_2=\nabla\cdot\mathbf{a}y+\mathbf{a}\cdot\nabla y
+\nabla B\cdot\mathbf{z}+B\na\cdot\mathbf{z}$.
 Similar to the estimate of $\mathbf{f_1}$,
 by (\ref{e4.3}), we get
 $f_2\in L^p(0,T;L^{p_1}(\O))$, and
 \bel{e4.5}
 \|f_2\|_{L^p(0,T;L^{p_1}(\O))}
 \leq CM_2\(\|(y,\mathbf{z})\|_{V^2\times V^2}
 +\|\mathbf{h}\|_{L^p(\o\times (0,T))}+
 \|\mathbf{z}_0\|_{W^{2-\frac{2}{p},p}(\O)}\).
 \ee
 Applying Theorem 2.3 in \cite{2006-18} again to the solution $y$, we deduce that
 $$y\in L^p(0,T;W^{2,p_1}(\O)), \ y_t\in L^p(0,T;L^{p_1}(\O))$$ and
 \begin{equation*}
 \|y\|_{L^p(0,T;W^{2,p_1}(\O))}
 +\|y_t\|_{L^p(0,T;L^{p_1}(\O))}
 \leq  C\(\|f_2\|_{L^p(0,T;L^{p_1}(\O))}
 +\|y_0\|_{W^{2-\frac{2}{p},p}(\O)}\).
 \end{equation*}
 By (\ref{e4.5}), it follows that
 \begin{eqnarray*}
 &&\|y\|_{L^p(0,T;W^{2,p_1}(\O))}
 +\|y_t\|_{L^p(0,T;L^{p_1}(\O))}\\
 &&\leq CM_2\(\|(y,\mathbf{z})\|_{V^2\times V^2}
 +\|y_0\|_{W^{2-\frac{2}{p},p}(\O)}
 +\|\mathbf{z}_0\|_{W^{2-\frac{2}{p},p}(\O)}
 +\|\mathbf{h}\|_{L^p(\o\times (0,T))}\).
 \end{eqnarray*}
 If $p\leq \frac{2np}{np-4}$, i.e., $p\leq 2+\frac{4}{n}$, this ends the proof.
 If $p>\frac{2np}{np-4}$, the proof will be completed by repeating the above procedure for finitely many times.
 \endpf

 The null controllability result for the equation (\ref{1.7}) can be stated as follows.

\begin{theorem}\label{t1}
 Assume that the condition (\ref{e82}) holds and $T>0$. Then there exists a function $\mathbf{h}\in L^2(Q)$ such that the associated solution $(y,\mathbf{z})$ of equations (\ref{1.7}) satisfies
$$y(x,T)=\mathbf{z}(x,T)=0, \ \text{a.e.}\ x\in\O.$$
Moreover,
\bel{1.8}
\|\mathbf{h}\|_{L^2(Q)}\leq
C\(\|y_0\|_{L^2(\O)}+\|\mathbf{z}_0\|_{L^2(\O)}\).
\ee
\end{theorem}

We consider the following adjoint system of (\ref{1.7}):
\begin{equation}\label{1.9}
\left\{\begin{array}{ll}
-\varphi_t-\Delta \varphi+\mathbf{a}\cdot\nabla\varphi
=-\nabla\cdot\bm{\psi},& (x,t)\in Q,\\
\bm{\psi}_t+\Delta\bm{\psi}
+\mathbf{b}\nabla\cdot\bm{\psi}
=B\nabla\varphi, &(x,t)\in Q,\\
\varphi=\bm{\psi}=0,&(x,t)\in\Sigma,\\
(\varphi,\bm{\psi})(x,T)=(\varphi_0,\bm{\psi}_0)(x),
& x\in\Omega.
\end{array}\right.
\end{equation}

 We derive a new global Carleman estimate for (\ref{1.9}).
  Assume that
 $\rho\in C^2(\overline{\O})$ satisfies
 \bel{1.10}
 |\nabla\rho|\geq C>0\  \text{in}\  \O\backslash\overline{\o_0}, \ \rho>0\
  \text{in} \ \O, \ \text{and}\  \rho=0 \ \text{on} \ \partial\O,
   \ee
   where $\o_0\neq\emptyset$ is an open subset of $\o$. Let $\o_1$ be any
   fixed open subset of $\o$ such that $\overline{\o_0}\subseteq \o_1$ and $\overline{\o_1}\subseteq \o$.
   Inspired by \cite{guerrero}, we first introduce the weight functions. For any real number $\l>1$ and $s>1$,
set
$$\theta=e^l,\ l=-s\phi,\ \phi(x,t)=\frac{\exp\{\frac{k(m+1)}{m}\l \|\rho\|_{L^\infty(\O)}\}-\exp\{\l (k\|\rho\|_{L^\infty(\O)}+\rho(x))\}}{t^m(T-t)^m},$$
where $m>3$ and $k>m$ are fixed.
 Put
 $$\xi(x,t)=\frac{\exp\{\l (k\|\rho\|_{L^\infty(\O)}+\rho(x))\}}{t^m(T-t)^m},$$
 $$\phi^*(t)=\max\limits_{x\in\overline{\O}}\phi(x,t)
 =\phi(x,t)\big{|}_{\partial\O},\
 \xi^*(t)=\min\limits_{x\in\overline{\O}}\xi(x,t)
 =\xi(x,t)\big{|}_{\partial\O}.$$

We have the following global Carleman estimate for (\ref{1.9}).
\begin{theorem}\label{t2}
Assume that the condition (\ref{e82}) holds. Then
there exist $\l_1,s_1>0$ such that for all $\l\geq \l_1,\ s\geq s_1$, one can find a constant $C>0$
 such that the following inequality holds for the solutions of (\ref{1.9}):
 \begin{eqnarray}\label{1.11}
 &&s\l^2\int_Q\theta^2\xi|\Delta\varphi|^2dxdt
 +s^3\l^4\int_Q\theta^2\xi^3|\nabla\varphi|^2dxdt
 +s^6\l^8\int_Q\theta^2\xi^6|\bm{\psi}|^2dxdt\nonumber\\
&&\leq C(1+T^{2m})s^8\l^8\int_0^T\int_\o\theta^4e^{2s\phi^*}\xi^8|\bm{\psi}|^2dxdt.
\end{eqnarray}
\end{theorem}

Before giving the proof of Theorem \ref{t2}, we first recall a Carleman estimate for the parabolic equation with nonhomogeneous  Neumann boundary conditions, which will be useful (see \cite{cara}).

\begin{lemma}\label{l2}
Let $\overline{y}_0\in L^2(\O)$, $h_1\in L^2(Q)$, $\bm{h_2}\in L^2(Q)$
 and $h_3\in L^2(\Sigma)$. Then there is a constant $C(\O,\o_0)>0$, such that for any $\l\geq C$ and $s\geq C(T^{2m}+T^{2m-1})$, any solution $\overline{y}\in L^2(0,T;H^1(\O))\cap L^\infty(0,T;L^2(\O))$ of
 \begin{equation*}
\left\{\begin{array}{ll}
\overline{y}_t-\Delta \overline{y}=h_1+\nabla\cdot \bm{h_2}, &(x,t)\in Q,\\
\partial_\nu \overline{y}+\bm{h_2}\cdot\nu=h_3,&(x,t)\in\Sigma,\\
\overline{y}(x,0)=\overline{y}_0(x),
& x\in\Omega
\end{array}\right.
\end{equation*}
satisfies that
\begin{eqnarray*}
&&s\l^2\int_Q\theta^2\xi|\nabla \overline{y}|^2dxdt
+s^3\l^4\int_Q\theta^2\xi^3|\overline{y}|^2dxdt\\
&&\leq C\(s^3\l^4\int_0^T\int_{\o_0}\theta^2\xi^3|\overline{y}|^2dxdt+
\int_Q\theta^2|h_1|^2dxdt
+s^2\l^2\int_Q\theta^2\xi^2|\bm{h_2}|^2dxdt\\[2mm]
&&\quad+s\l\int_\Sigma e^{-2s\phi^*}\xi^*|h_3|^2d\sigma dt\).
\end{eqnarray*}
\end{lemma}

\no{\bf Proof of Theorem \ref{t2}.}
The main idea of this proof is borrowed from \cite{guerrero}.
The proof will be divided into two steps.
\medskip

\no{\bf Step 1.} We first consider the parabolic equation satisfied by $\nabla\varphi$, because $B\nabla\varphi$ appears on the right hand side of the  equation satisfied by $\bm{\psi}$.
 By (\ref{1.9}), we know that $\nabla \varphi$ satisfies
 \bel{1.12}
 -(\nabla \varphi)_t-\Delta (\nabla \varphi)+\nabla(\mathbf{a}\cdot\nabla \varphi)=
 -\nabla(\nabla\cdot\bm{\psi})\ \text{in}\ Q.
 \ee
 Set $\mathbf{a}=(a_1,a_2,...,a_n)^{\top}$, $\bm{\psi}=(\psi_1,\psi_2,...,\psi_n)^{\top}$,
then by (\ref{1.12}), it follows that $\frac{\partial\varphi}{\partial x_i}$ satisfies
$$-\(\frac{\partial\varphi}{\partial x_i}\)_t-\Delta\(\frac{\partial\varphi}{\partial x_i}\)
+\sum\limits_{j=1}^n\(a_j\frac{\partial^2\varphi}{\partial x_i\partial x_j}+\frac{\partial a_j}{\partial x_i}\frac{\partial\varphi}{\partial x_j}\)
=-\sum\limits_{j=1}^n
  \frac{\partial^2\psi_j}{\partial x_i\partial x_j}
  =-\nabla\cdot \frac{\partial\bm{\psi}}{\partial x_i}.$$
   Applying Lemma \ref{l2} for $\frac{\partial\varphi}{\partial x_i}$, here
  $h_1=-\sum\limits_{j=1}^n\(a_j\frac{\partial^2\varphi}{\partial x_i\partial x_j}+\frac{\partial a_j}{\partial x_i}\frac{\partial\varphi}{\partial x_j}\)$, $\bm{h_2}=-\frac{\partial\bm{\psi}}{\partial x_i}$ and
$h_3=\sum\limits_{j=1}^n\frac{\partial^2\varphi}{\partial x_i\partial x_j}\nu_j-\sum\limits_{j=1}^n\frac{\partial\psi_j}{\partial \nu}\nu_i\nu_j$, we conclude that

\begin{eqnarray}\label{1.14}
&&I(\nabla\varphi):=s\l^2\int_Q\theta^2\xi\sum\limits_{i=1}^n\Big{|}\nabla \(\frac{\partial\varphi}{\partial x_i}\)\Big{|}^2dxdt
+s^3\l^4\int_Q\theta^2\xi^3\sum\limits_{i=1}^n\Big{|}\frac{\partial\varphi}{\partial x_i}\Big{|}^2dxdt\nonumber\\
&&\ \ \ \ \ \ \ \ \ \leq C\[s^3\l^4\displaystyle\int_0^T\!\!\int_{\o_0}
\theta^2\xi^3\sum\limits_{i=1}^n\Big{|}\frac{\partial\varphi}{\partial x_i}\Big{|}^2dxdt+s^2\l^2\int_Q\theta^2\xi^2
\sum\limits_{i,j=1}^n\Big{|}\frac{\partial\psi_j}{\partial x_i}\Big{|}^2dxdt
\nonumber\\
&&\quad\quad\quad\quad+s\l\int_\Sigma e^{-2s\phi^*}\xi^*\sum\limits_{i=1}^n\Big{|}\sum\limits_{j=1}^n\(\frac{\partial^2\varphi}{\partial x_i\partial x_j}\nu_j-\frac{\partial\psi_j}{\partial \nu}\nu_i\nu_j\)\Big{|}^2d\sigma dt\nonumber\\
&&\quad\quad\quad\quad+\int_Q\theta^2\sum\limits_{i=1}^n
\Big{|}\sum\limits_{j=1}^n\(a_j\frac{\partial^2\varphi}{\partial x_i\partial x_j}+\frac{\partial a_j}{\partial x_i}\frac{\partial\varphi}{\partial x_j}\)\Big{|}^2dxdt\].
\end{eqnarray}

Notice that $\frac{1}{\xi}\leq CT^{2m}$, and by the definition of $\xi$, it follows that
\begin{eqnarray}\label{1.15}
&&\int_Q\theta^2\sum\limits_{i=1}^n
\Big{|}\sum\limits_{j=1}^n\(a_j\frac{\partial^2\varphi}{\partial x_i\partial x_j}+\frac{\partial a_j}{\partial x_i}\frac{\partial\varphi}{\partial x_j}\)\Big{|}^2dxdt\nonumber\\
&&\leq C \(\|\mathbf{a}\|^2_{L^{\infty}(Q)}+ \|\na\mathbf{a}\|^2_{L^{\infty}(Q)}\)
\[T^{2m} \int_Q\theta^2\xi\sum\limits_{i=1}^n\Big{|} \na\(\frac{\partial\varphi}{\partial x_i}\)\Big{|}^2dxdt \nonumber\\
&&\ \
+ T^{6m}\int_Q \theta^2\xi^3\sum\limits_{i=1}^n\Big{|}
\frac{\partial\varphi}{\partial x_i}\Big{|}^2dxdt\].
\end{eqnarray}
We next estimate $s\l\int_\Sigma e^{-2s\phi^*}\xi^*\sum\limits_{i=1}^n\Big{|}\sum\limits_{j=1}^n\(\frac{\partial^2\varphi}{\partial x_i\partial x_j}\nu_j-\frac{\partial\psi_j}{\partial \nu}\nu_i\nu_j\)\Big{|}^2d\sigma dt$. To do this, take $$\rho(t)=s^{\frac{1}{2}-\frac{1}{m}}\l e^{-s\phi^*(t)}(\xi^*)^{\frac{1}{2}-\frac{1}{m}}(t)$$ and define
$\varphi^*=\rho(t)\varphi$, then $\varphi^*$ will be the solution of
 \begin{equation}\label{1.16}
\left\{\begin{array}{ll}
-\varphi_t^*-\Delta \varphi^*+\mathbf{a}\cdot\nabla\varphi^*=-\rho\nabla\cdot\bm{\psi}
-\rho_t\varphi, &(x,t)\in Q,\\
\varphi^*=0,&(x,t)\in\Sigma,\\
\varphi^*(x,T)=0,
& x\in\Omega.
\end{array}\right.
\end{equation}
It is easy to check that $-\rho\nabla\cdot\bm{\psi}
-\rho_t\varphi\in L^2(0,T;H^1(\O))$, then $\varphi^*\in
L^2(0,T;H^3(\O))\cap H^1(0,T;H^1(\O))$, and
\bea\label{1.17}
\|\varphi^*\|^2_{L^2(0,T;H^3(\O))}
+\|\varphi^*_t\|^2_{L^2(0,T;H^1(\O))}\leq C\(\|\rho_t\varphi\|^2_{L^2(0,T;H^1(\O))}
+\|\rho\nabla\cdot\bm{\psi}\|^2_{L^2(0,T;H^1(\O))}\).
\eea
Moreover,
\bea\label{1.18}
&&\|\rho_t\varphi\|^2_{L^2(0,T;H^1(\O))}
\leq C\int_Q\rho_t^2\sum\limits_{i=1}^n\Big{|}\frac{\partial\varphi}{\partial x_i}\Big{|}^2dxdt
\leq
CTs^{\frac{3}{2}-\frac{1}{m}}\l\int_Qe^{-2s\phi^*}(\xi^*)^3
\sum\limits_{i=1}^n\Big{|}\frac{\partial\varphi}{\partial x_i}\Big{|}^2dxdt\nonumber\\
&&
\leq Cs^{\frac{3}{2}}\l \int_Qe^{-2s\phi^*}(\xi^*)^3
\sum\limits_{i=1}^n\Big{|}\frac{\partial\varphi}{\partial x_i}\Big{|}^2dxdt
\leq CI(\nabla\varphi)
\eea
for $s\geq CT^m$.
By (\ref{1.17}) and (\ref{1.18}), we get
\bea\label{1.19}
 &&\|\varphi^*\|^2_{L^2(0,T;H^3(\O))}=s^{1-\frac{2}{m}}\l^2\int_0^T
e^{-2s\phi^*}(\xi^*)^{1-\frac{2}{m}}\|\varphi\|^2_{H^3(\O)}dt\nonumber\\[2mm]
 &&\leq C\(I(\nabla\varphi)+
\|\rho\nabla\cdot\bm{\psi}\|^2_{L^2(0,T;H^1(\O))}\).
\eea
 From this, using the integration by parts, we conclude that
\begin{eqnarray}\label{h1.20}
s^{2-\frac{1}{m}}\l^3\int_0^T
e^{-2s\phi^*}(\xi^*)^{2-\frac{1}{m}}\|\varphi\|^2_{H^2(\O)}dt
\leq C\(I(\nabla\varphi)+
\|\rho\nabla\cdot\bm{\psi}\|^2_{L^2(0,T;H^1(\O))}\).
\end{eqnarray}
Combining (\ref{1.19}) with (\ref{h1.20}) yields
\begin{eqnarray}\label{1.21}
s^{\frac{3}{2}-\frac{3}{2m}}\l^3\int_0^T
e^{-2s\phi^*}(\xi^*)^{\frac{3}{2}-\frac{3}{2m}}
\sum\limits_{i,j=1}^n\Big{\|}\frac{\partial^2\varphi}{\partial x_i\partial x_j}\nu_j\Big{\|}^2_{L^2(\Sigma)}dt\leq C\(I(\nabla\varphi)+
\|\rho\nabla\cdot\bm{\psi}\|^2_{L^2(0,T;H^1(\O))}\).
\end{eqnarray}

In addition,  by the trace theorem,
\bel{1.22}
s\l\int_\Sigma e^{-2s\phi^*}\xi^*\sum\limits_{i,j=1}^n\Big{|}\frac{\partial\psi_j}{\partial \nu}\nu_i\nu_j\Big{|}^2d\sigma dt
\leq Cs\l \int_0^T e^{-2s\phi^*}\xi^*\|\bm{\psi}\|^2_{H^2(\O)}dt.
\ee
By (\ref{1.14}), (\ref{1.15}), (\ref{1.21}) and (\ref{1.22}), noting that $m>3$, we deduce that
there exists a constant $C>0$ such that, for any $s\geq\max\{CT^{2m},CT^m\}$ and $\l\geq C \(\|\mathbf{a}\|_{L^{\infty}(Q)} + \|\na\mathbf{a}\|_{L^{\infty}(Q)}\)$, it holds that
\begin{eqnarray}\label{1.23}
&&I(\nabla\varphi)
\leq C\[s^3\l^4\displaystyle\int_0^T\!\!\int_{\o_0}
\theta^2\xi^3\sum\limits_{i=1}^n\Big{|}\frac{\partial\varphi}{\partial x_i}\Big{|}^2dxdt+s^2\l^2\int_Q\theta^2\xi^2
\sum\limits_{i,j=1}^n\Big{|}\frac{\partial\psi_j}{\partial x_i}\Big{|}^2dxdt
\nonumber\\
&&\quad\quad\quad\quad+s\l^2\int_0^T e^{-2s\phi^*}\xi^*\|\nabla\bm{\psi}\|^2_{H^1(\O)} dt
\].
\end{eqnarray}

{\bf Step 2.} Let us consider the equation satisfied by $\bm{\psi}$. Write $\bm{b}=(b_1,b_2,...,b_n)^{\top}$
and
$$J(\bm{\psi}):=s^6\l^8\int_Q\theta^2\xi^6|\bm{\psi}|^2dxdt
+s^4\l^6\int_Q\theta^2\xi^4|\nabla\bm{\psi}|^2dxdt
+s^2\l^4\int_Q\theta^2\xi^2|\Delta\bm{\psi}|^2dxdt.$$
By (\ref{1.9}), applying the classical Carleman estimate of the parabolic operator with the right-hand side in $L^2(Q)$ for $\bm{\psi}$, we see that
\bea\label{1.24}
J(\bm{\psi})\leq
C\(s^6\l^8\int_0^T\int_{\o_0}\theta^2\xi^6|\bm{\psi}|^2dxdt
+s^3\l^4\|B\|_{L^\infty(Q)}^2\int_Q\theta^2\xi^3|\nabla\varphi|^2dxdt\).
\eea
Multiplying (\ref{1.23}) by $\(1+\|B\|_{L^\infty(Q)}^2\)$, and adding (\ref{1.24}) to it, we conclude that, for any $s\geq\max\{CT^{2m},CT^m\}$  and $\l\geq C\(\|\mathbf{a}\|_{L^{\infty}(Q)}+ \|\na\mathbf{a}\|_{L^{\infty}(Q)}+\|B\|_{L^{\infty}(Q)}\)$, it holds that
\bea\label{1.25}
&&I(\nabla\varphi)+J(\bm{\psi})
\leq C\[s^3\l^4\int_0^T\int_{\o_0}\theta^2\xi^3|\nabla\varphi|^2dxdt
+s\l^2\int_0^Te^{-2s\phi^*}\xi^*\|\nabla\bm{\psi}\|^2_{H^1(\O)}dt\nonumber\\
&&\quad\quad\quad\quad\quad\quad\quad+s^6\l^8\int_0^T\int_{\o_0}\theta^2\xi^6|\bm{\psi}|^2dxdt\].
\eea
We next claim that the second term in the right hand side of (\ref{1.25}) can be absorbed by the left hand side.
To this end, we set $\zeta(t)=s^{\frac{1}{2}}\l e^{-s\phi^*}(\xi^*)^{\frac{1}{2}}$ and $\bm{\Psi}=\zeta(t)\bm{\psi}$.
Then $\bm{\Psi}$ satisfies
\begin{equation*}
\left\{\begin{array}{ll}
\bm{\Psi}_t+\Delta \bm{\Psi}+\bm{b}\nabla\cdot\bm{\Psi}=\zeta B\nabla\varphi+\zeta_t\bm{\psi}, &(x,t)\in Q,\\
\bm{\Psi}=0, &(x,t)\in\Sigma,\\
\bm{\Psi}(x,T)=0,
& x\in\Omega.
\end{array}\right.
\end{equation*}
By a simple calculation, we have $|\zeta_t(t)|\leq Ts^{\frac{3}{2}}\l e^{-s\phi^*}(\xi^*)^{\frac{3}{2}+\frac{2}{m}}$ and  $\zeta B\nabla\varphi+\zeta_t\bm{\psi}\in L^2(Q)$, then $\bm{\Psi}\in L^2(0,T;H^2(\O))$ and
\bea\label{1.26}
&&\|\bm{\Psi}\|^2_{L^2(0,T;H^2(\O))}
=s\l^2\int_0^Te^{-2s\phi^*}\xi^*\|\bm{\psi}\|^2_{H^2(\O)}dt
\nonumber\\
&&\leq C\(s\l^2\|B\|^2_{L^\infty(Q)}\int_Qe^{-2s\phi^*}\xi^*|\nabla\varphi|^2dxdt
+T^2s^3\l^2\int_Qe^{-2s\phi^*}(\xi^*)^{3+\frac{4}{m}}|\bm{\psi}|^2dxdt\)
\nonumber\\
&&\leq
C\(I(\nabla\varphi)+J(\bm{\psi})\)
\eea
for any $\l\geq C\|B\|_{L^\infty(Q)}$.
By (\ref{1.25}) and (\ref{1.26}), we deduce that
\bea\label{1.27}
 &&I(\nabla\varphi)+J(\bm{\psi})
+\|\bm{\Psi}\|^2_{L^2(0,T;H^2(\O))}\nonumber\\[2mm]
 &&\leq C\(s^3\l^4\int_0^T\int_{\o_0}\theta^2\xi^3|\nabla\varphi|^2dxdt
+s^6\l^8\int_0^T\int_{\o_0}\theta^2\xi^6|\bm{\psi}|^2dxdt\)
\eea
for any $s\geq\max\{CT^{2m},CT^m\}$  and $\l\geq C(\|\mathbf{a}\|_{L^{\infty}(Q)}+ \|\na\mathbf{a}\|_{L^{\infty}(Q)}+\|B\|_{L^\infty(Q)})$.

\medskip
 We proceed to show that the first term in the right hand side of (\ref{1.27}) can be also eliminated.
Notice that $$\bm{\psi}_t+\Delta\bm{\psi}
+\mathbf{b}\nabla\cdot\bm{\psi}
=B\nabla\varphi\ \text{in}\ \o_0\times(0,T).$$
Take $\varrho\in C_0^2(\o)$ and $\varrho\equiv1$ in $\o_0$, where $w_0\subset\subset \o$.
Then, integrating by parts, we obtain
 \bea\label{1.28}
 &&s^3\l^4\int_0^T\int_{\o}\varrho\theta^2\xi^3|B\nabla\varphi|^2dxdt
 =s^3\l^4\int_0^T\int_{\o}\varrho\theta^2\xi^3B\nabla\varphi\(\bm{\psi}_t+\Delta\bm{\psi}
+\mathbf{b}\nabla\cdot\bm{\psi}\)dxdt
\nonumber\\
&&=-s^3\l^4\int_0^T\int_{\o}\varrho(\theta^2\xi^3)_t
B\nabla\varphi\cdot\bm{\psi}dxdt
+ s^3\l^4\int_0^T\int_{\o}\Delta(\varrho\theta^2\xi^3)B\nabla\varphi\cdot\bm{\psi}dxdt
\nonumber\\
&&+2 s^3\l^4\int_0^T\int_{\o}\nabla(\varrho\theta^2\xi^3)\nabla(B\nabla\varphi)\cdot\bm{\psi}dxdt
- s^3\l^4\int_0^T\int_{\o}\nabla(\varrho\theta^2\xi^3)
B\nabla\varphi\bm{b}\cdot\bm{\psi}dxdt
\nonumber\\
&&+ s^3\l^4\int_0^T\int_{\o}\varrho\theta^2\xi^3\bm{\psi}
\(
\Delta(B\nabla\varphi)-(B\nabla\varphi)_t-
\nabla(B\nabla\varphi\cdot\bm{b})\)dxdt.
\eea
An easy verification shows that
$$(\theta^2\xi^3)_t\leq CTs\theta^2\xi^{4+\frac{1}{m}},\
\nabla(\theta^2\xi^3)\leq Cs\l\theta^2\xi^4\
\text{and}\ \Delta(\theta^2\xi^3)\leq Cs^2\l^2\theta^2\xi^5.$$
Applying these, and Young's inequality, we get
\bea\label{1.29}
&&-s^3\l^4\int_0^T\int_{\o}\varrho(\theta^2\xi^3)_t
B\nabla\varphi\cdot\bm{\psi}dxdt
\leq CTs^4\l^4\int_0^T\int_{\o}\theta^2\xi^{4+\frac{1}{m}}
\varrho B\nabla\varphi\cdot\bm{\psi}dxdt\nonumber\\
&&\leq \varepsilon_0C(\|B\|^2_{L^\infty(Q)})I(\nabla\varphi)
+C(\varepsilon_0)s^5\l^4T^2\int_0^T\int_{\o}
\varrho^2\theta^2
\xi^{5+\frac{2}{m}}
|\bm{\psi}|^2dxdt,
\eea
\bea\label{1.30}
&& s^3\l^4\int_0^T\int_{\o}
\Delta(\varrho\theta^2\xi^3)B\nabla\varphi
\cdot\bm{\psi}dxdt
\leq Cs^5\l^6\int_0^T\int_{\o}
\varrho\theta^2\xi^5B\nabla\varphi
\cdot\bm{\psi}dxdt
\nonumber\\
&&\leq \varepsilon_0C(\|B\|^2_{L^\infty(Q)})I(\nabla\varphi)
+C(\varepsilon_0)s^7\l^8\int_0^T\int_{\o}
\varrho^2\theta^2
\xi^{7}
|\bm{\psi}|^2dxdt,
\eea
\bea\label{1.31}
&&2 s^3\l^4\int_0^T\int_{\o}\nabla(\varrho\theta^2\xi^3)\nabla(B\nabla\varphi)\cdot\bm{\psi}dxdt
\leq Cs^4\l^5\int_0^T\int_{\o}
\varrho\theta^2\xi^4|\nabla B\cdot\nabla\varphi+B\na(\na\varphi)||\bm{\psi}|dxdt
\nonumber\\
&&\leq
\varepsilon_0C(\|B\|^2_{L^\infty(Q)}+\|\na B\|^2_{L^\infty(Q)})I(\nabla\varphi)
+C(\varepsilon_0)s^7\l^8\int_0^T\int_{\o}
\varrho^2\theta^2
\xi^{7}
|\bm{\psi}|^2dxdt,
\eea
and
\bea\label{1.32}
&&-s^3\l^4\int_0^T\int_{\o}\nabla(\varrho\theta^2\xi^3)
B\nabla\varphi\bm{b}\cdot\bm{\psi}dxdt
\leq Cs^4\l^5\int_0^T\int_{\o}|\varrho\theta^2\xi^4
B\nabla\varphi\bm{b}\cdot\bm{\psi}|dxdt
\nonumber\\
&&\leq
\varepsilon_0C(\|B\|^2_{L^\infty(Q)}+\|\bm{b}\|^2_{L^\infty(Q)})I(\nabla\varphi)
+C(\varepsilon_0)s^5\l^6\int_0^T\int_{\o}
\varrho^2\theta^2
\xi^{5}
|\bm{\psi}|^2dxdt.
\eea
By (\ref{1.9}), we have
\bea\label{h2}
-(\nabla\varphi)_t=\Delta(\nabla\varphi)-\nabla\mathbf{a}\cdot\nabla\varphi
-\mathbf{a}\cdot\nabla(\nabla\varphi)-\nabla(\nabla\cdot\bm{\psi}).
\eea
Then,  we have
\bea\label{1.33}
&&\Delta(B\nabla\varphi)-(B\nabla\varphi)_t-
\nabla(B\nabla\varphi\cdot\bm{b})\nonumber\\
&&=(\Delta  B-B_t- \nabla B\cdot\bm{b}- B\nabla\bm{b}-B\nabla\mathbf{a})\nabla\varphi
\nonumber\\
&&\quad+(2\nabla B-B\mathbf{a}- B\bm{b})\cdot\nabla(\nabla\varphi)
+2B\Delta(\nabla\varphi)-B\nabla(\nabla\cdot\bm{\psi}).
\eea
From this, we get
\bea\label{1.34}
&&s^3\l^4\int_0^T\int_{\o}\varrho\theta^2\xi^3\bm{\psi}
\(
\Delta(B\nabla\varphi)-(B\nabla\varphi)_t-
\nabla(B\nabla\varphi\cdot\bm{b})\)dxdt
\nonumber\\
&&=s^3\l^4\int_0^T\int_{\o}\varrho\theta^2\xi^3\bm{\psi}
(\Delta  B-B_t- \nabla B\cdot\bm{b}- B\nabla\bm{b}-B\nabla\mathbf{a})\nabla\varphi dxdt\nonumber\\
&&\quad
+s^3\l^4\int_0^T\int_{\o}\varrho\theta^2\xi^3\bm{\psi}
(2\nabla B-B\mathbf{a}-B\bm{b})\cdot\nabla(\nabla\varphi)dxdt
\nonumber\\
&&\quad+2s^3\l^4\int_0^T\int_{\o}\varrho\theta^2\xi^3\bm{\psi}
B\Delta(\nabla\varphi)dxdt
-s^3\l^4\int_0^T\int_{\o}\varrho\theta^2\xi^3\bm{\psi}
B\nabla(\nabla\cdot\bm{\psi})dxdt\nonumber\\
&&:=I_1+I_2+I_3+I_4.
\eea

In what follows,
let $C_0$ denote a constant dependent on $\|B\|_{L^\infty(Q)}$,  $\|B_t\|_{L^\infty(Q)}$, $\|\na B\|_{L^\infty(Q)}$,\ $\|\Delta B\|_{L^\infty(Q)}$,\ $\|\bm{b}\|_{L^\infty(Q)}$,
 $\|\na\bm{b}\|_{L^\infty(Q)}$, $\|\mathbf{a}\|_{L^\infty(Q)}$ and $\|\na\mathbf{a}\|_{L^\infty(Q)}$, which may vary from line to line.
By Young's inequality, we obtain
\bea\label{1.35}
I_1\leq
\varepsilon_0C_0I(\nabla\varphi)+C(\varepsilon_0)
s^3\l^4\int_0^T\int_{\o}\varrho^2\theta^2\xi^3|\bm{\psi}|^2dxdt,
\eea
\bea\label{1.36}
I_2\leq
\varepsilon_0C_0I(\nabla\varphi)+C(\varepsilon_0)
s^5\l^6\int_0^T\int_{\o}\varrho^2\theta^2\xi^5|\bm{\psi}|^2dxdt,
\eea
\bea\label{1.37}
&&I_3=-2s^3\l^4\int_0^T\int_{\o}\(
\nabla(\varrho\theta^2\xi^3)\cdot\bm{\psi}
B\Delta\varphi +\varrho\theta^2\xi^3\nabla\cdot\bm{\psi}B\Delta\varphi
+\varrho\theta^2\xi^3\bm{\psi}\cdot\nabla B\Delta\varphi\) dxdt
\nonumber\\
&&\quad\leq
C(\varepsilon_0)
s^7\l^8\int_0^T\int_{\o}\varrho^2\theta^2\xi^7|\bm{\psi}|^2dxdt
+C(\varepsilon_0)
s^5\l^6\int_0^T\int_{\o}
\varrho^2\theta^2\xi^5|\nabla\bm{\psi}|^2dxdt\nonumber\\[2mm]
&&\quad\quad+\varepsilon_0C_0I(\nabla\varphi),
\eea
\bea\label{1.38}
&&I_4
\leq \varepsilon_0C_0s\l^2\int_0^T\int_{\o}e^{-2s\phi^*}\xi^*
|\nabla(\nabla\cdot \bm{\psi})|^2dxdt
+C(\varepsilon_0)
s^5\l^6\int_0^T\int_{\o}
\varrho^2\theta^4e^{2s\phi^*}
\frac{\xi^6}{\xi^*}|\bm{\psi}|^2dxdt
\nonumber\\
&&\quad\leq \varepsilon_0C_0\|\bm{\Psi}\|^2_{L^2(0,T;H^2(\O))}
+C(\varepsilon_0)T^{2m}
s^5\l^6\int_0^T\int_{\o}
\varrho^2e^{-4s\phi+2s\phi^*}\xi^6|\bm{\psi}|^2dxdt.
\eea
By (\ref{1.34})--(\ref{1.38}), we conclude that
\bea\label{1.39}
&&s^3\l^4\int_0^T\int_{\o}\varrho\theta^2\xi^3\bm{\psi}
\(
\Delta(B\nabla\varphi)-(B\nabla\varphi)_t-
\nabla(B\nabla\varphi\cdot\bm{b})\)dxdt
\nonumber\\
&&\leq \varepsilon_0C_0\(I(\nabla\varphi)+\|\bm{\Psi}\|^2_{L^2(0,T;H^2(\O))}\)
+C(\varepsilon_0)T^{2m}
s^7\l^8\int_0^T\int_{\o}
\varrho^2e^{-4s\phi+2s\phi^*}\xi^7|\bm{\psi}|^2dxdt
\nonumber\\
&&\quad
+C(\varepsilon_0)
s^5\l^6\int_0^T\int_{\o}
\varrho^2\theta^2\xi^5|\nabla\bm{\psi}|^2dxdt.
\eea
Moreover, we have
\bea\label{1.40}
&&s^5\l^6\int_0^T\int_{\o}
\varrho^2\theta^2\xi^5|\nabla\bm{\psi}|^2dxdt
=s^5\l^6\int_0^T\int_{\o}\(
\nabla(\varrho^2\theta^2\xi^5)\nabla\bm{\psi}\bm{\psi}+
\varrho^2\theta^2\xi^5\nabla(\nabla\bm{\psi})\bm{\psi}\)dxdt
\nonumber\\
&&\leq
\varepsilon_0J(\bm{\psi})+C(\varepsilon_0)s^8\l^8
\int_0^T\int_{\o}
\varrho^4\theta^2\xi^8|\bm{\psi}|^2dxdt.
\eea
Combining (\ref{1.39}) with (\ref{1.40}) yields
\bea\label{1.41}
&&s^3\l^4\int_0^T\int_{\o}\varrho\theta^2\xi^3\bm{\psi}
\(
\Delta(B\nabla\varphi)-(B\nabla\varphi)_t-
\nabla(B\nabla\varphi\cdot\bm{b})\)dxdt
\nonumber\\
&&
\leq
\varepsilon_0C_0\(I(\nabla\varphi)+J(\bm{\psi})+\|\bm{\Psi}\|^2_{L^2(0,T;H^2(\O))}\)
\nonumber\\
&&\quad+C(\varepsilon_0)T^{2m}
s^8\l^8\int_0^T\int_{\o}
\varrho^2e^{-4s\phi+2s\phi^*}\xi^8|\bm{\psi}|^2dxdt.
\eea
Note that $B$ has a positive lower bound, by (\ref{1.28})-(\ref{1.32}) and (\ref{1.41}), we obtain
\bea\label{1.42}
&&s^3\l^4\int_0^T\int_{\o_0}
\theta^2\xi^3|\nabla\varphi|^2dxdt
\leq Cs^3\l^4\int_0^T\int_{\o_0}\varrho
\theta^2\xi^3|B\nabla\varphi|^2dxdt\nonumber\\
&&\leq Cs^3\l^4\int_0^T\int_{\o}\varrho
\theta^2\xi^3|B\nabla\varphi|^2dxdt
\leq
\varepsilon_0C_0\(I(\nabla\varphi)+J(\bm{\psi})+\|\bm{\Psi}\|^2_{L^2(0,T;H^2(\O))}\)
\nonumber\\
&&\quad+C(\varepsilon_0)(1+T^{2m})
s^8\l^8\int_0^T\int_{\o}
\varrho^2e^{-4s\phi+2s\phi^*}\xi^8|\bm{\psi}|^2dxdt.
\eea

Substituting (\ref{1.42}) into (\ref{1.27}), and choosing $\varepsilon_0$ small enough, we conclude that
\bea\label{1.43}
I(\nabla\varphi)+J(\bm{\psi})
+\|\bm{\Psi}\|^2_{L^2(0,T;H^2(\O))}\leq C(1+T^{2m})s^8\l^8\int_0^T\int_{\o}
e^{-4s\phi+2s\phi^*}\xi^8|\bm{\psi}|^2dxdt,
\eea
which completes the proof.
\endpf

By the classical fact, the statement of Theorem \ref{t1}
will be obtained once we prove the lemma below, which can be obtained by Theorem \ref{t2}.

\begin{lemma}\label{l1}
Assume that the condition (\ref{e82}) holds. Then
there exists a constant $C=C(\O,\o,T)>0$ independent of
$(\varphi_0,\bm{\psi}_0)$ such that
\bel{ee3}
\|\varphi|_{t=0}\|^2_{L^2(\O)}+\|\bm{\psi}|_{t=0}\|^2_{L^2(\O)}
\leq C\|\theta^2e^{s\phi^*}\xi^4\bm{\psi}\|^2_{L^2(\o\times(0,T))},
\ \forall (\varphi_0,\bm{\psi}_0)\in L^2(\O)\times L^2(\O).
\ee
\end{lemma}

\no{\bf Proof.}  Multiplying (\ref{1.9}) by $\varphi$ and $\bm{\psi}$ respectively, integrating over $\O
$, and using integration by parts, we have
\begin{eqnarray*}
\frac{1}{2}\frac{d}{dt}\int_\O(\varphi^2+\bm{\psi}^2)dx-\int_\O(|\na \varphi|^2+|\na\bm{\psi}|^2)dx
=\int_\O(B\na \varphi \bm{\psi}+\na\cdot \bm{\psi}\varphi-\mathbf{b}\nabla\cdot\bm{\psi}\bm{\psi}+\mathbf{a}\cdot\nabla\varphi\varphi)dx.
\end{eqnarray*}
By Cauchy inequality, we have
\begin{eqnarray*}
&&\frac{1}{2}\frac{d}{dt}\int_\O(\varphi^2+\bm{\psi}^2)dx-\int_\O(|\na \varphi|^2+|\na\bm{\psi}|^2)dx\\[2mm]
&&\geq-\frac{C}{\varepsilon}(1+\|B\|^2_{L^\infty(\O)}+\|\mathbf{b}\|^2_{L^\infty(\O)}
+\|\mathbf{a}\|^2_{L^\infty(\O)})\int_\O(\varphi^2+\bm{\psi}^2)dx
-\varepsilon \int_\O(|\na \varphi|^2+|\na\bm{\psi}|^2)dx.
\end{eqnarray*}
Taking $\varepsilon$ small enough, we obtain
\begin{eqnarray*}
\frac{1}{2}\frac{d}{dt}\int_\O(\varphi^2+\bm{\psi}^2)dx
\geq -C(1+\|B\|^2_{L^\infty(\O)}+\|\mathbf{b}\|^2_{L^\infty(\O)}
+\|\mathbf{a}\|^2_{L^\infty(\O)})\int_\O(\varphi^2+\bm{\psi}^2)dx.
\end{eqnarray*}
By Gronwall's inequality, we get
$$\int_\O(\varphi^2(x,0)+\bm{\psi}^2(x,0)dx\leq C\int_\O(\varphi^2(x,t)+\bm{\psi}^2(x,t))dx.$$
Integrating this inequality on $(\frac{T}{4},\frac{3T}{4})$, we have
\begin{eqnarray}\label{ee1}
\|\varphi|_{t=0}\|^2_{L^2(\O)}+\|\bm{\psi}|_{t=0}\|^2_{L^2(\O)}
\leq C\int^{\frac{3T}{4}}_{\frac{T}{4}}\int_\O(\varphi^2+\bm{\psi}^2)dxdt.
\end{eqnarray}

Moreover, by Poincar\'{e}'s inequality and the definitions of $\theta$ and $\xi$, we have
\begin{eqnarray}\label{ee2}
&&s^3\l^4\int_Q\theta^2\xi^3|\nabla\varphi|^2dxdt
 +s^6\l^8\int_Q\theta^2\xi^6|\bm{\psi}|^2dxdt\nonumber\\[2mm]
 &&\geq Cs^3\l^4\int^{\frac{3T}{4}}_{\frac{T}{4}}\int_\O\varphi^2dxdt+
 Cs^6\l^8\int^{\frac{3T}{4}}_{\frac{T}{4}}\int_\O\bm{\psi}^2dxdt.
\end{eqnarray}
By (\ref{ee1}), (\ref{ee2}) and Theorem \ref{t2}, we deduce (\ref{ee3}). This completes the proof of Lemma
\ref{l1}.
\endpf
%


\begin{remark}
Note that it is a technical condition that $B$ has a positive lower bound in Theorem \ref{t2}, which plays a critical role in  the proof of Theorem \ref{t2}. Indeed, in order to establish the global Carleman estimate for (\ref{1.9}), we need to derive a local estimate for $\na\varphi$ (see the first term in the right hand side of (\ref{1.27})). However, this local estimate for $\na\varphi$ is obtained by estimating $B\na\varphi$  because $B\na\varphi$  appears in the equation satisfied by $\bm{\psi}$.
Therefore, we require that $B$ has a positive lower bound, then (\ref{1.42}) and (\ref{1.43}) hold.
\end{remark}

\begin{remark}\label{r2}
Notice that in the proof of Theorem \ref{t2},  we require that the coefficients in system (\ref{1.9}) satisfy
$\mathbf{a}, \mathbf{b},$
$\nabla\mathbf{a},$ $\nabla\mathbf{b}\in  L^\infty(Q),$
 $B, B_t, \na B, \Delta B\in L^\infty(Q)$, by the relationship between $\mathbf{a}, \mathbf{b}, B$ and $\overline{\mathbf{v}}, \overline{u}$ (see (\ref{e09})) and the embedding theorem, we need to establish the regularity of $(\ou,\ov)$, which actually lies in $C([0,T];H^4(\O))$
with smallness on $H^3$ norm of the initial data.
\end{remark}

%

Next, we establish the null controllability of system  (\ref{1.7}) with a control function in $L^r(Q)$, where $r>n+2$.

%

\begin{proposition}\label{p2}
Assume that the condition (\ref{e82}) holds. Let $r>n+2$, $y_0,\mathbf{z}_0\in W^{2-\frac{2}{r},r}(\O)\cap H_0^1(\O)$. Then one can find a control $\mathbf{h}\in L^r(Q)$ supported in $\o\times [0,T]$
  such that the solution $(y,\mathbf{z})\in V^r\times V^r$  of system (\ref{1.7}) satisfies
  \bel{6.2}
  y(x,T)=\mathbf{z}(x,T)=0\ \text{a.e.\ in}\ \O.
  \ee
 Moreover,
\bel{6.1}
\|\mathbf{h}\|_{L^r(Q)}\leq C e^{CM_1}M_2^{k_1}(\|y_0\|_{L^2(\O)}
+\|\mathbf{z}_0\|_{L^2(\O)}).
\ee
\end{proposition}

We give the following lemma, which will be needed in the proof of Proposition \ref{p2}.
Let $m_0,\gamma\geq 1$. Consider the following Banach space
\begin{equation}\label{e4.16}
X^{m_0,\gamma}(Q):=
L^\infty(0,T;L^{m_0}(\O))\cap L^\gamma(0,T;W^{1,\gamma}(\O)),
\end{equation}
 equipped with the norm
$\|v\|_{X^{m_0,\gamma}(Q)}=\esssup\limits_{0<t<T}
\|v(\cdot,t)\|_{L^{m_0}(\O)}+\|Dv\|_{L^\gamma(Q)}$.


\begin{lemma}\label{l5}(\cite[Proposition 3.2]{dege})
There exists a constant $C>0$ depending only upon $n,\gamma$ and $m_0$ such that for every $v\in X^{m_0,\gamma}(Q)$, it holds that
$$\|v\|_{L^q(Q)}\leq C\(1+\frac{T}{|\O|^{\frac{n(\gamma-m_0)+m_0\gamma}{nm_0}}}\)^{\frac{1}{q}}
\|v\|_{X^{m_0,\gamma}(Q)},$$
where $q=\gamma \frac{n+m_0}{n}.$
\end{lemma}

\no{\bf Proof of Proposition \ref{p2}.} For any given $\epsilon>0$, we consider the following optimal control problem:
\bel{e4.6}
(P_\epsilon): \ \text{Min}\{\frac{1}{2}\int_Qe^{4s\phi-2s\phi^*}
\xi^{-8}|\mathbf{h}|^2dxdt+\frac{1}{2\epsilon}
\int_{\O}y^2(x,T)dx+\frac{1}{2\epsilon}
\int_{\O}|\mathbf{z}(x,T)|^2dx, \ \text{subject to}\ (\ref{1.7})
\}.
\ee
By the standard variational method, we know that
for any $\epsilon>0$, the problem $(P_\epsilon)$ has a unique solution $(y_\epsilon, \mathbf{z}_\epsilon,\mathbf{h}_\epsilon)$ and
\bel{e4.8}
\mathbf{h}_\epsilon=\chi_\o e^{-4s\phi+2s\phi^*}
\xi^{8}\bm{\psi}_\epsilon,
\ee
where $(\varphi_\epsilon,\bm{\psi}_\epsilon)$ satisfies
\begin{equation}\label{e4.7}
\left\{\begin{array}{ll}
-\varphi_{\epsilon,t}-\Delta \varphi_\epsilon+\mathbf{a}\cdot\nabla\varphi_\epsilon
=-\nabla\cdot\bm{\psi}_\epsilon,& (x,t)\in Q,\\
\bm{\psi}_{\epsilon,t}+\Delta\bm{\psi}_\epsilon
+\mathbf{b}\nabla\cdot\bm{\psi}_\epsilon
=B\nabla\varphi_\epsilon, &(x,t)\in Q,\\
\varphi_\epsilon=\bm{\psi}_\epsilon=0,&(x,t)\in\Sigma,\\
(\varphi_\epsilon,\bm{\psi}_\epsilon)(x,T)=
-\frac{1}{\epsilon}(y_\epsilon(x ,T),\mathbf{z}_\epsilon(x,T)),
& x\in\Omega,
\end{array}\right.
\end{equation}
where $(y_\epsilon,\mathbf{z}_\epsilon)$ is the solution of (\ref{1.7}) associated to $(y_0,\mathbf{z}_0)$ and $\mathbf{h}_\epsilon$.

Multiplying the first (resp. second) equation of (\ref{e4.7})
 by $y_\epsilon$ (resp.  $\mathbf{z}_\epsilon)$ and integrating it on $Q$, by the boundary conditions of (\ref{1.7}) and (\ref{e4.7}), we have
 \begin{eqnarray*}
 &&\int_{\O}y_\epsilon(x,T)\varphi_\epsilon(x,T)dx
 -\int_{\O}y_\epsilon(x,0)\varphi_\epsilon(x,0)dx
 +\int_{\O}\bm{\psi}_\epsilon(x,T)\mathbf{z}_\epsilon(x,T)dx
 -\int_{\O}\bm{\psi}_\epsilon(x,0)\mathbf{z}_\epsilon(x,0)dx
 \\
 &&=\int_Q\chi_\o\mathbf{h}_\epsilon\bm{\psi}
 _\epsilon dxdt.
 \end{eqnarray*}
From (\ref{e4.8}), (\ref{e4.7}) and Lemma \ref{l1}, we obtain
 \begin{eqnarray*}
&&\frac{1}{\epsilon}\int_\O y_\epsilon^2(x,T)dx
+\frac{1}{\epsilon}\int_\O \mathbf{z}_\epsilon^2(x,T)dx
+\int_Q \chi_\o^2e^{-4s\phi+2s\phi^*}\xi^8
\bm{\psi}_\epsilon^2dxdt\\
&&=-\int_\O y_0(x)\varphi_\epsilon(x,0)dx
-\int_\O \mathbf{z}_0(x)\bm{\psi}_\epsilon(x,0)dx\\
&&\leq C\(\|y_0\|_{L^2(\O)}+\|\mathbf{z}_0\|_{L^2(\O)}\)\cdot
\left\|\bm{\psi}_\epsilon e^{-2s\phi+s\phi^*}\xi^4\right\|_{L^2(\o\times(0,T))},
\end{eqnarray*}
which implies that
\bea\label{e4.10}
&&\frac{1}{\epsilon}\int_\O y_\epsilon^2(x,T)dx
+\frac{1}{\epsilon}\int_\O \mathbf{z}_\epsilon^2(x,T)dx
+\int_Q \chi_\o e^{-4s\phi+2s\phi^*}\xi^8
\bm{\psi}_\epsilon^2dxdt\nonumber\\[2mm]
&&\leq C\(\|y_0\|^2_{L^2(\O)}+\|\mathbf{z}_0\|^2_{L^2(\O)}\),
\eea
which implies that $\{\mathbf{h}_\epsilon\}$ is a family of ``approximate" control, because  $y_\epsilon(x,T)\rightarrow 0$ in $L^2(\O)$, $\mathbf{z}_\epsilon(x,T)\rightarrow 0$ in $L^2(\O)$
as $\epsilon\rightarrow 0$.

Next, we prove that $\mathbf{h}_\epsilon\in L^r(Q)$ with $r>n+2$.
Let $\t>0$ and $\{\t_k\}_{k\in\mathbb{N}}$ be an increasing sequence such that
$0<\t_k<\t<\frac{s}{2}$.
Set
\bel{e4.15}
\bm{\Phi}_\epsilon^k=e^{-(s+\t_k)\phi^*(t)}
(\xi^*)^8(t)\bm{\psi}_\epsilon,\
\Psi_\epsilon^k=e^{-(s+\t_k)\phi^*(t)}
(\xi^*)^8(t)\varphi_\epsilon.\
\ee
On one hand, it is easy to prove that $\bm{\Phi}_\epsilon^k$
satisfies
\begin{equation}\label{e4.11}
\left\{\begin{array}{ll}
\bm{\Phi}_{\epsilon,t}^k+\Delta \bm{\Phi}_\epsilon^k+
\mathbf{b}\na\cdot\bm{\Phi}_\epsilon^k\\[1.5mm]=
\[e^{-(s+\t_k)\phi^*}(\xi^*)^8\]_t\bm{\psi}_\epsilon
+\[e^{-(s+\t_k)\phi^*}(\xi^*)^8\]B\na\varphi_\epsilon
:=\mathbf{g}_\epsilon^k,& (x,t)\in Q,\\
\bm{\Phi}_\epsilon^k=0,&(x,t)\in\Sigma,\\
\bm{\Phi}_\epsilon^k(x,0)=\bm{\Phi}_\epsilon^k(x,T)=0,
& x\in\Omega.
\end{array}\right.
\end{equation}
Then
\bea\label{e4.12}
\|\mathbf{g}_\epsilon^1\|^2_{L^2(Q)}
&\leq& 2\int_Q\left|\[e^{-(s+\t_1)\phi^*}(\xi^*)^8\]_t\right|^2\bm{\psi}_\epsilon^2dxdt
+2\int_Q\[e^{-(s+\t_1)\phi^*}(\xi^*)^8\]^2B^2|\na\varphi_\epsilon|^2dxdt
\nonumber\\&:=&I_1+I_2.
\eea
Since
\begin{eqnarray*}
&&\Big{|}\[e^{-(s+\t_1)\phi^*}(\xi^*)^8\]_t\Big{|}
=\Big{|}-(s+\t_1)\p_t^*e^{-(s+\t_1)\p^*}(\xi^*)^8
+8e^{-(s+\t_1)\p^*}(\xi^*)^7\xi^*_t\Big{|}\\
&&\leq (s+\t_1)C(T)(\xi^*)^{9+\frac{1}{m}}e^{-(s+\t_1)\p^*}
+C(T)e^{-(s+\t_1)\p^*}(\xi^*)^{8+\frac{1}{m}},
\end{eqnarray*}
we obtain
\begin{eqnarray*}
&&I_1\leq C\int_Qs^2(\xi^*)^{2(9+\frac{1}{m})}e^{-2(s+\t_1)\p^*}\bm{\psi}_\epsilon^2dxdt
+C\int_Qe^{-2(s+\t_1)\p^*}(\xi^*)^{2(8+\frac{1}{m})}\bm{\psi}_\epsilon^2dxdt
\\&&\leq Cs^2\int_Qe^{-2s\p}\xi^6\bm{\psi}_\epsilon^2dxdt,
\end{eqnarray*}
since $(\xi^*)^{12+\frac{2}{m}}e^{-2\t_1\p^*}\leq1$ and
 $(\xi^*)^{10+\frac{2}{m}}e^{-2\t_1\p^*}\leq1$.

Similarly, we have
\begin{eqnarray*}
I_2\leq C(\|B\|_{L^\infty(Q)})
\int_Qe^{-2s\p}\xi^3|\na\varphi_\epsilon|^2dxdt.
\end{eqnarray*}
Therefore, by Theorem \ref{t2} and (\ref{e4.10}), it follows that
$$\|\mathbf{g}_\epsilon^1\|^2_{L^2(Q)}\leq
 C
\int_0^T\int_{\o}e^{-4s\p+2s\p^*}\xi^8
\bm{\psi}_\epsilon^2dxdt
\leq C\(\|y_0\|_{L^2(\O)}^2+\|\mathbf{z}_0\|_{L^2(\O)}^2\).
$$
By Proposition \ref{p1}, we have $\bm{\Phi}_\epsilon^1\in V^2$. Moreover,
$$\|\bm{\Phi}_\epsilon^1\|_{V^2}^2
\leq e^{CM_1}M_2^{k_1}\|\mathbf{g}_\epsilon^1\|^2_{L^2(Q)}
\leq e^{CM_1}M_2^{k_1}\(\|y_0\|_{L^2(\O)}^2+\|\mathbf{z}_0\|_{L^2(\O)}^2\).$$
By the embedding theorem, $V^2=W_2^{2,1}(Q)\hookrightarrow L^{s_1}(Q)$ for
 $s_1=\left\{\begin{array}{ll}
\frac{2(n+2)}{n-2},& n>2,\\
\text{any constant}\  \kappa>1,&n\leq 2.
\end{array}\right.$ Then,
$$\|\bm{\Phi}_\epsilon^1\|_{L^{s_1}(Q)}^2
\leq e^{CM_1}M_2^{k_1}
\(\|y_0\|_{L^2(\O)}^2+\|\mathbf{z}_0\|^2_{L^2(\O)}\).$$

On the other hand, it is easy to check that
$\Psi_\epsilon^k$ satisfies
\begin{equation}\label{e4.12}
\left\{\begin{array}{ll}
-{\Psi}_{\epsilon,t}^k-\Delta \Psi_\epsilon^k+\mathbf{a}\cdot
\na\Psi_\epsilon^k\\[1.5mm]=
-\[e^{-(s+\t_k)\phi^*}(\xi^*)^8\]_t\varphi_\epsilon
-\[e^{-(s+\t_k)\phi^*}(\xi^*)^8\]\na\cdot\bm{\psi}_\epsilon
:=f_\epsilon^k,& (x,t)\in Q,\\
\Psi_\epsilon^k=0,&(x,t)\in\Sigma,\\
\Psi_\epsilon^k(x,0)=\Psi_\epsilon^k(x,T)=0,
& x\in\Omega.
\end{array}\right.
\end{equation}
Next, we prove $f_\epsilon^1\in L^2(Q)$. Similarly, by (\ref{e4.10}), Theorem \ref{t2} (or (\ref{1.43})) and Poincar\'{e}'s inequality,  we deduce
\begin{eqnarray*}
&&\int_Q|f_\epsilon^1|^2dxdt
\leq Cs^3\l^4\int_Qe^{-2s\p}\xi^3|\na\varphi_\epsilon|^2dxdt
+s^4\l^6\int_Qe^{-2s\p}\xi^4|\na\cdot\bm{\psi}_\epsilon|^2dxdt
\\&&\leq Cs^8\l^8\int_0^T\int_{\o}
e^{-4s\p+2s\p^*}\xi^8|\bm{\psi}_\epsilon|^2dxdt
\leq C\(\|y_0\|_{L^2(\O)}^2+\|\mathbf{z}_0\|_{L^2(\O)}^2\).
\end{eqnarray*}
By Proposition \ref{p1}, we know $\Psi_\epsilon^1\in V^2$. Moreover,
$$\|\Psi_\epsilon^1\|_{V^2}^2
\leq e^{CM_1}M_2^{k_1}\|f_\epsilon^1\|^2_{L^2(Q)}
\leq  e^{CM_1}M_2^{k_1}
\(\|y_0\|_{L^2(\O)}^2+\|\mathbf{z}_0\|_{L^2(\O)}^2\).$$
 By the embedding theorem, $V^2\hookrightarrow L^{s_1}(Q)$, then
 $$\|\Psi_\epsilon^1\|_{L^{s_1}(Q)}^2
\leq e^{CM_1}M_2^{k_1}
\(\|y_0\|_{L^2(\O)}^2+\|\mathbf{z}_0\|_{L^2(\O)}^2\).$$

In what follows, we give the estimates of $\mathbf{g}_\epsilon^2$ and $f_\epsilon^2$, respectively.
By (\ref{e4.15}) and (\ref{e4.11}), we have
\bea\label{e4.13}
&&\mathbf{g}_\epsilon^2=\[e^{-(s+\t_2)\phi^*}(\xi^*)^8\]_t\bm{\psi}_\epsilon
+\[e^{-(s+\t_2)\phi^*}(\xi^*)^8\]B\na\varphi_\epsilon
\nonumber\\&&\quad\ =\[e^{-(s+\t_2)\phi^*}(\xi^*)^8\]_t
e^{(s+\t_1)\phi^*}(\xi^*)^{-8} \bm{\Phi}_\epsilon^1
+Be^{(\t_1-\t_2)\p^*}\na\Psi_\epsilon^1.
\eea

Notice that, by (\ref{e4.12}) and Proposition \ref{p1}, we have
$\na \Psi_{\epsilon}^1\in L^\infty(0,T;L^2(\O))\cap L^2(0,T; W^{1,2}(\O))$. Taking $m_0=2, \gamma=2$ in (\ref{e4.16}), by Lemma \ref{l5}, we deduce that
$\na \Psi_{\epsilon}^1\in L^{q_1}(Q)$, where $q_1=\frac{2(n+2)}{n}>2$, and
\bel{e4.18}
\|\na \Psi_{\epsilon}^1\|_{L^{q_1}(Q)}^2
\leq C\|\na \Psi_{\epsilon}^1\|_{X^{2,2}(Q)}^2
\leq Ce^{CM_1}M_2^{k_1}
\(\|y_0\|_{L^2(\O)}^2+\|\mathbf{z}_0\|_{L^2(\O)}^2\).
\ee
Moreover,
 \bel{e4.17}
 \Big{|}\[e^{-(s+\t_2)\phi^*}(\xi^*)^8\]_t
e^{(s+\t_1)\phi^*}(\xi^*)^{-8}\Big{|}
\leq Ce^{(\tau_1-\tau_2)}\phi^*(\xi^*)^{\frac{1}{m}}
\ \text{and}\ \bm{\Phi}_\epsilon^1\in L^{s_1}(Q),
 \ee
 where we choose $s_1=q_1=\frac{2(n+2)}{n}$ when $n=1$ or $2$, and $s_1=\frac{2(n+2)}{n-2}>q_1$ when $n=3$.
Therefore,
By (\ref{e4.13})-(\ref{e4.17}), we get that
\bel{e4.117}
\|\mathbf{g}_\epsilon^2\|_{L^{q_1}(Q)}^2
\leq Ce^{CM_1}M_2^{k_1}
\(\|y_0\|_{L^2(\O)}^2+\|\mathbf{z}_0\|_{L^2(\O)}^2\).
\ee
Again, by (\ref{e4.11}) and Proposition \ref{p1}, we see that
$\bm{\Phi}_\epsilon^2\in V^{q_1}$.
Moreover,
$$\|\bm{\Phi}_\epsilon^2\|_{V^{q_1}}^2
\leq Ce^{CM_1}M_2^{k_1}\|\mathbf{g}_\epsilon^2\|_{L^{q_1}(Q)}^2
\leq Ce^{CM_1}M_2^{k_1}
\(\|y_0\|_{L^2(\O)}^2+\|\mathbf{z}_0\|_{L^2(\O)}^2\).$$
By the embedding theorem, $V^{q_1}\hookrightarrow L^{s_2}(Q)$, where $s_2=\left\{\begin{array}{ll}
\frac{q_1(n+2)}{n+2-2q_1},& n+2-2q_1>0,\\
\text{any constant}\  \kappa>1,&n+2-2q_1\leq 0.
\end{array}\right.$
Hence,
\bel{e4.19}
\|\bm{\Phi}_\epsilon^2\|_{L^{s_2}(Q)}^2
\leq Ce^{CM_1}M_2^{k_1}
\(\|y_0\|_{L^2(\O)}^2+\|\mathbf{z}_0\|_{L^2(\O)}^2\).
\ee
In addition, by (\ref{e4.15}) and (\ref{e4.12}), we arrive at
$$f_\epsilon^2=-[e^{-(s+\tau_2)\phi^*}(\xi^*)^8]_t
e^{(s+\tau_1)\phi^*}(\xi^*)^{-8}\Psi_\epsilon^1
-e^{(\tau_1-\tau_2)\phi^*}
\na\cdot\bm\Phi_\epsilon^1.$$
Similar to (\ref{e4.18}), we can prove that $\na\cdot\bm\Phi_\epsilon^1\in L^{q_1}(Q)$. Combining with
$\Psi_\epsilon^1\in L^{s_1}(Q)$, we have $f_\epsilon^2\in L^{q_1}(Q).$
 Using Proposition \ref{p1} again, we deduce that
 $\Psi_\epsilon^2\in V^{q_1}$.

 By the embedding theorem, it follows that
 \bel{e4.21}
 \|\Psi_\epsilon^2\|_{L^{s_{2}}(Q)}^2
 \leq C\|\Psi_\epsilon^2\|_{V^{q_1}}^2
 \leq Ce^{CM_1}M_2^{k_1}
\(\|y_0\|_{L^2(\O)}^2+\|\mathbf{z}_0\|_{L^2(\O)}^2\).
 \ee
 Similarly, since $\na \Psi_\epsilon^2\in L^\infty(0,T;L^2(\O))\cap L^{q_1}(0,T;W^{1,q_1}(\O))$,
 by Lemma \ref{l5}, we take $m_0=2,\ \gamma=q_1$, it follows that
 \bel{e4.20}
 \na \Psi_\epsilon^2\in L^{q_2}(Q), \ \text{where}\
 q_2=\frac{q_1(n+2)}{n}.
 \ee
 Combining with $\bm\Phi_\epsilon^2\in L^{s_2}(Q)$, we have $\mathbf{g}_\epsilon^3\in L^{q_2}(Q)$ and
 $$\|\mathbf{g}_\epsilon^3\|_{L^{q_2}(Q)}^2
 \leq Ce^{CM_1}M_2^{k_1}
\(\|y_0\|_{L^2(\O)}^2+\|\mathbf{z}_0\|_{L^2(\O)}^2\),$$
where, we take $s_2=q_2$, since $n+2-2q_1<0$.
 By (\ref{e4.11}), it follows that
 $$\bm\Phi_\epsilon^3\in V^{q_2}\hookrightarrow
  L^{s_3}(Q),
 \ \text{where}\  s_3=\left\{\begin{array}{ll}
\frac{q_2(n+2)}{n+2-2q_2},& n+2-2q_2>0,\\
\text{any constant}\  \kappa>1&n+2-2q_2\leq 0.
\end{array}\right.$$

 Moreover, similar to (\ref{e4.20}), we can show that
 $\na\cdot \bm\Phi_\epsilon^2\in L^{q_2}(Q)$. By (\ref{e4.21}), $\Psi_\epsilon^2\in L^{s_2}(Q)$. Then,
 $f_\epsilon^3\in L^{q_2}(Q)$. By (\ref{e4.12}),
  $\Psi_\epsilon^3\in V^{q_2}\hookrightarrow
  L^{s_3}(Q)$.

  Repeating the above procedure, since $q_{N+1}-q_N=
  q_N\(\frac{n+2}{n}-1\)=q_N\cdot \frac{2}{n}>0$,
   there exists a $N^*\in \mathbb{N}$ such that
  $$ \bm\Phi_\epsilon^{N^*}\in L^{q_{N^*}}(Q),\ \Psi_\epsilon^{N^*}\in L^{q_{N^*}}(Q),\ \text{where}\
  q_{N^*}>n+2.$$
  By (\ref{e4.8}),
  $$\mathbf{h}_\epsilon=\chi_\o e^{-4s\phi+2s\phi^*}
\xi^{8}\bm{\psi}_\epsilon
=\chi_\o e^{-4s\phi+2s\phi^*}
\xi^{8}e^{(s+\tau_{N^*})\phi^*}
(\xi^*)^{-8}\bm\Phi_\epsilon^{N^*}.$$
 Since $\tau_{N^*}<\tau<\frac{s}{2},$ one has
  $e^{-4s\phi+2s\phi^*}
\xi^{8}e^{(s+\tau_{N^*})\phi^*}
(\xi^*)^{-8}\leq C$. Hence,
$\mathbf{h}_\epsilon \in L^r(Q),$ where $r>n+2$.
 Moreover,
 \bel{e4.22}
 \|\mathbf{h}_\epsilon\|_{L^r(Q)}
 \leq Ce^{CM_1}M_2^{k_1}
\(\|y_0\|_{L^2(\O)}+\|\mathbf{z}_0\|_{L^2(\O)}\).
 \ee
 Letting $\epsilon\rightarrow 0$, by (\ref{e4.22}) and (\ref{e4.10}), we conclude that
  there exists a control $\mathbf{h}\in L^r(Q)$
  such that
  the solution of (\ref{1.7}) satisfies
  $y(x,T)=\mathbf{z}(x,T)=0$ in $\O$.
  Moreover,
$$
 \|\mathbf{h}\|_{L^r(Q)}
 \leq Ce^{CM_1}M_2^{k_1}
\(\|y_0\|_{L^2(\O)}+\|\mathbf{z}_0\|_{L^2(\O)}\),
 $$
 which is the desired conclusion.\endpf

 \section{The Proof of main result}

 \no{\bf Proof of Theorem \ref{t4}.}
 Set
 $K=\{\bm{\eta}\in V^r\ |\ \|\bm{\eta}\|_{V^r}\leq 1\}.$
 For any $\bm{\eta}\in K$, we consider the following linearized system:
 \begin{equation}\label{e09}
\left\{\begin{array}{ll}
y_t-\Delta y=\nabla\cdot(\mathbf{a_{\bm{\eta}}}y)
+\nabla\cdot(B\mathbf{z}),& (x,t)\in Q,\\
\mathbf{z}_t-\Delta\mathbf{z}=
-\nabla(\mathbf{b_{\bm{\eta}}}\cdot\mathbf{z})
+\nabla y+
\chi_{\o}\mathbf{h}, &(x,t)\in Q,\\
y=\mathbf{z}=0,&(x,t)\in\Sigma,\\
(y,\mathbf{z})(x,0)=(y_0,\mathbf{z}_0)(x),
& x\in\Omega,
\end{array}\right.
\end{equation}
where $\mathbf{a_{\bm{\eta}}}=\bm{\eta}+\overline{\mathbf{v}}$,
$\mathbf{b_{\bm{\eta}}}=\bm{\eta}+2\overline{\mathbf{v}}$ and $B=\overline{u}$.
  Define
 \begin{eqnarray*}
 &&\Lambda(\bm{\eta})=\{\mathbf{z}\in V^r|\ \exists\  \mathbf{h}\in L^r(Q)\  \text{and a constant}\ C>0\  \text{such that the solution of}\ (\ref{e09})\  \\
 &&\ \ \ \ \ \ \ \ \ \ \  \text{corresponding to} \ \bm{\eta} \ \text{and}\ \mathbf{h}\  \text{satisfies}\
 (\ref{6.2})\ \text{and}\ (\ref{6.1})\}.
 \end{eqnarray*}
 Obviously, $K$ is a nonempty convex subset of $V^r$.
 By Proposition \ref{p2}, we know that $\Lambda(\bm{\eta})$  is a nonempty convex subset of  $V^r$.

 Next, we prove that $\Lambda(\bm{\eta})$  is a compact subset of $V^r$.  By Proposition \ref{p1} and (\ref{6.1}), we get
 \bel{6.3}
 \|\mathbf{z}\|_{V^r}
 \leq e^{CM_1}M_2^{k_1}\(\|(y_0,\mathbf{z}_0)\|_{W^{2-\frac{2}{r},r}(\O)
 \times W^{2-\frac{2}{r},r}(\O)}+\|y_0\|_{L^2(\O)}
 +\|\mathbf{z}_0\|_{L^2(\O)}\).
 \ee
 Therefore, $\|\mathbf{z}\|_{V^r}$ is bounded.
  Note that, when $r>n+2$, $V^r\hookrightarrow C^{1+\alpha,\frac{1+\alpha}{2}}(Q),$ here $\alpha=1-\frac{n+2}{r}$.
  Applying the Arzela-Ascoll Theorem, we can obtain that
 $\Lambda(\bm{\eta})$  is a compact subset of $V^r$.

 Further, we show that $\Lambda$ is upper semi-continuous.
 For this, let $\{\bm{\eta}_n\}_{n=1}^\infty\subset K$ such that $\bm{\eta}_n\rightarrow \bm{\eta}$ in $K$, and set
 $\mathbf{z}_n\in \Lambda(\bm{\eta}_n)$.
 By the definition of $\Lambda(\bm{\eta}_n)$, there exists
 $\mathbf{h}_n\in L^r(Q)$ such that the solution $(y_n,\mathbf{z}_n)$ of (\ref{e09}) satisfies (\ref{6.2}) and (\ref{6.1}).
 By Proposition \ref{p1}, we  have
 $$ \|y_n\|_{V^r}+\|\mathbf{z}_n\|_{V^r}
 \leq e^{CM_1}M_2^{k_1}\(\|(y_0,\mathbf{z}_0)\|_{W^{2-\frac{2}{r},r}(\O)
 \times W^{2-\frac{2}{r},r}(\O)}+\|y_{0}\|_{L^2(\O)}
 +\|\mathbf{z}_{0}\|_{L^2(\O)}\).$$
 Hence, there exist $\mathbf{h}\in L^r(Q)$, $y,\mathbf{z}\in V^r$, and the subsequences of
 $\{\mathbf{h}_n\}, \{y_n\}, \{\mathbf{z}_n\}$ (still denoted by themselves), such that
 \bel{6.6}
 \mathbf{h}_n\rightharpoonup \mathbf{h}\ \text{in}\ L^r(Q),
 \ y_n\rightharpoonup y\ \text{in}\ V^r,
 \ \text{and}\
 \mathbf{z}_n\rightharpoonup \mathbf{z}\ \text{in}\ V^r.
 \ee
 Then $(y,\mathbf{z})$ is the solution of (\ref{e09}) corresponding to $\bm{\eta}$ and $\mathbf{h}$.
 Take $Y_n=y_n-y$, $\mathbf{Z}_n=\mathbf{z}_n-\mathbf{z}$,
 and $\mathbf{H}_n=\chi_\o(\mathbf{h}_n-\mathbf{h})$. Then $(Y_n,\mathbf{Z}_n)$ satisfies
 \begin{equation}\label{6.4}
\left\{\begin{array}{ll}
Y_{n,t}-\Delta Y_n=\na\cdot[\mathbf{a}_{\bm{\eta_n}}Y_n+(\mathbf{a}_{\bm{\eta_n}}
-\mathbf{a}_{\bm{\eta}})y]
+\na\cdot(B\mathbf{Z}_n),& (x,t)\in Q,\\
\mathbf{Z}_{n,t}-\Delta \mathbf{Z}_n=
-\na[\mathbf{b}_{\bm{\eta_n}}\cdot\mathbf{Z}_n+
(\mathbf{b}_{\bm{\eta_n}}
-\mathbf{b}_{\bm{\eta}})\cdot\mathbf{z}]
+\na Y_n +\mathbf{H}_n, &(x,t)\in Q,\\
Y_n=\mathbf{Z}_n=0,&(x,t)\in\Sigma,\\
(Y_n,\mathbf{Z}_n)(x,0)=(0,0),
& x\in\Omega.
\end{array}\right.
\end{equation}
 Moreover, an easy computation shows that
 \begin{eqnarray}\label{6.5}
 &&\|Y_n(\cdot,t)\|^2_{L^2(\O)}+
  \|\na Y_n(\cdot,t)\|^2_{L^2(\O)}+
   \|\mathbf{Z}_n(\cdot,t)\|^2_{L^2(\O)}\\\nonumber
 &&\leq e^{M_2T}\(\int_\O \mathbf{H}_n\cdot \mathbf{Z}_ndx
 +\int_\O|\bm{\eta}_n-\bm{\eta}|^2(|y|^2+|\mathbf{z}|^2)dx
 \).
\end{eqnarray}
By (\ref{6.6}), it follows that
\bel{6.7}
\mathbf{Z}_n\rightarrow 0 \ \text{in}\ L^{r_0}(Q),\ \text{where}\ r_0=\frac{r}{r-1}.
\ee
 By (\ref{6.5}) and (\ref{6.7}), we have
 $$ \|Y_n(\cdot,t)\|^2_{L^2(\O)}\rightarrow 0,\ \|\mathbf{Z}_n(\cdot,t)\|^2_{L^2(\O)}\rightarrow 0,\ \forall\ t\in[0,T],$$
 and notice that $y_n(x,T)=\mathbf{z}_n(x,T)=0$ in $\O$.
 Hence, $y(x,T)=\mathbf{z}(x,T)=0$ in $\O$. i.e., $\mathbf{z}\in \Lambda(\bm{\eta})$.

 At last,  we claim that $\Lambda(\bm{\eta})\subset K$. Indeed,
 \begin{eqnarray*}
 &&\|y\|_{L^\infty(Q)}+\|\mathbf{z}\|_{V^r}
 \leq e^{CM_1}M_2^{k_1}\(\|(y_0,\mathbf{z}_0)\|_{W^{2-\frac{2}{r},r}(\O)
 \times W^{2-\frac{2}{r},r}(\O)}+\|y_0\|_{L^2(\O)}
 +\|\mathbf{z}_0\|_{L^2(\O)}\)\\
 &&\leq C\|(y_0,\mathbf{z}_0)\|_{W^{2-\frac{2}{r},r}(\O)
 \times W^{2-\frac{2}{r},r}(\O)}.
 \end{eqnarray*}
 Therefore, there exists a constant $\delta>0$ such that,
  if $\|(y_0,\mathbf{z}_0)\|_{W^{2-\frac{2}{r},r}(\O)
 \times W^{2-\frac{2}{r},r}(\O)}<\delta$, we have
 \bel{6.9}
 \|y\|_{L^\infty(Q)}+\|\mathbf{z}\|_{V^r(Q)}
 \leq C\|(y_0,\mathbf{z}_0)\|_{W^{2-\frac{2}{r},r}(\O)
 \times W^{2-\frac{2}{r},r}(\O)}
 \leq \min\{1, \overline{p}\}.
 \ee
 Thus, by the Kakutani's fixed point theorem, there exists
 $\mathbf{z}\in K$ such that $\mathbf{z}\in \Lambda(\mathbf{z})$.
 Moreover, by Corollary \ref{c22}, $\overline{u}\geq \overline{p}>0$. Therefore, $u=y+\overline{u}\geq 0$,
 which proves  Theorem \ref{t4}.
 \endpf

\section*{Acknowledgement}
The authors would like to thank the referees for valuable comments and suggestions. Tao is partially supported by the National Science Foundation of China under
grant 11971320 and Guangdong Basic and Applied Basic Research Foundation under
grant 2020A1515010530.
Zhang is partially supported by
the National Science Foundation of China under
grants 12001094, 12001087 and 11971179, and Fundamental Research Funds for the Central Universities under grant 2412020QD027.

\section*{Declarations}
\textbf{Conflict of interest} The authors certify that they have no interest directly or indirectly related to the work
submitted for publication.


\begin{thebibliography}{1}{\small}


\bibitem{AL1987} W. Alt, D.A. Lauffenburger, {\it Transient behavior of a chemotaxis system modelling certain types of tissue inflammation,} J. Math. Biol., 24(1987), 691--722.






\bibitem{2009427} F. Ammar-Khodja, A. Benabdallah, C. Dupaix, M. Gonz\'{a}lez-Burgos, {\it A generalization of the Kalman rank condition for
time-dependent coupled linear parabolic systems}, Differ. Equ. Appl., 1(2009), 427--457.

\bibitem{2009267} F. Ammar-Khodja, A. Benabdallah, C. Dupaix,  M. Gonz\'{a}lez-Burgos, {\it A Kalman
rank condition for the localized distributed controllability of a class of linear parbolic systems},
J. Evol. Equ., 9(2009), 267--291.

\bibitem{2005426} F. Ammar-Khodja, A. Benabdallah, C. Dupaix, I. Kostin, {\it Null-controllability of some systems of parabolic type by one
control force}, ESAIM Control Optim. Calc. Var.,  11(2005), 426--448.

\bibitem{FAML} F. Ammar-Khodja,  A. Benabdallah, M. Gonz\'{a}lez-Burgos,  L. de Teresa, {\it
Recent results on the controllability of linear coupled parabolic problems: a survey},
Math. Control Relat. Fields, 1(2011), 267--306.

 \bibitem{new}  F. Ammar-Khodja, A.  Benabdallah, M. Gonz\'{a}lez-Burgos, L. de Teresa, {\it
New phenomena for the null controllability of parabolic systems: minimal time and geometrical dependence,}
J. Math. Anal. Appl., 444(2016), 1071--1113.




\bibitem{BM1985} D. Balding, D.L. McElwain, {\it A mathematical model of tumour-induced capillary growth,} J. Theoret. Biol., 114(1985), 53--73.

\bibitem{barbu2} V. Barbu, {\it
Controllability of parabolic and Navier-Stokes equations},
Sci. Math. Jpn., 56(2002), 143--211.

\bibitem{yi}  A. Benabdallah, M. Cristofol, P. Gaitan, L. de Teresa, {\it
A new Carleman inequality for parabolic systems with a single observation and applications},
C. R. Math. Acad. Sci. Paris, 348(2010), 25--29.

\bibitem{201417} A. Benabdallah, M. Cristofol, P. Gaitan, L. de Teresa, {\it Controllability to trajectories for some parabolic systems of three
and two equations by one control force}, Math. Control Relat. Fields, 4(2014), 17--44.





\bibitem{CG2015} F.W. Chaves-Silva,  S. Guerrero, {\it
A uniform controllability result for the Keller-Segel system},
Asymptot. Anal., 92(2015), 313--338.

\bibitem{CG2017} F.W. Chaves-Silva, S. Guerrero, {\it A controllability result for a chemotaxis-fluid model}, J.
Differential Equations, 262(2017), 4863--4905.

\bibitem{CMS2020} J.M. Coron, F. Marbach, F. Sueur, {\it Small-time global exact controllability of the Navier-Stokes equation with Navier slip-with-friction boundary conditions}, J. Eur. Math. Soc., 22(2020), 1625--1673.

\bibitem{DLK1972} F.W. Dahlquist, P. Lovely, D.E. Koshland Jr, {\it Qualitative analysis of bacterial migration in chemotaxis,} Nat., New Biol., 236(1972), 120--123.


\bibitem{dege} E. DiBenedetto, {\it Degenerate parabolic equations}, Springer, New York, 2012.

\bibitem{si} M. Duprez, {\it
Controllability of a $2\times 2$ parabolic system by one force with space-dependent coupling term of order one},
ESAIM Control Optim. Calc. Var., 23(2017), 1473--1498.

\bibitem{ba} M. Duprez, P. Lissy, {\it Indirect controllability of some linear parabolic systems of m equations with m-1 controls involving
coupling terms of zero or first order}, J. Math. Pures Appl., 106(2016), 905--934.


\bibitem{shisan}  M. Duprez,  P. Lissy, {\it
Positive and negative results on the internal controllability of parabolic equations coupled by zero- and first-order terms},
J. Evol. Equ., 18(2018), 659--680.









 \bibitem{cara} E. Fern\'{a}ndez-Cara, M. Gonz\'{a}lez-Burgos, S. Guerrero, J.P. Puel, {\it
Null controllability of the heat equation with boundary Fourier conditions: the linear case},
ESAIM Control Optim. Calc. Var., 12(2006), 442--465.


\bibitem{2006-18} Y. Giga, H. Sohr, {\it
Abstract $L^p$ estimates for the Cauchy problem with applications to the Navier-Stokes equations in exterior domains},
J. Funct. Anal., 102(1991), 72--94.

\bibitem{201091} M. Gonz\'{a}lez-Burgos, L. de Teresa, {\it Controllability results for cascade systems of m coupled parabolic PDEs by one control
force}, Port. Math., 67(2010), 91--113.

\bibitem{guerrero} S. Guerrero, {\it
Null controllability of some systems of two parabolic equations with one control force},
SIAM J. Control Optim., 46(2007), 379--394.




\bibitem{GZ2014}  B.Z. Guo, L. Zhang, {\it
Local null controllability for a chemotaxis system of parabolic-elliptic type},
Systems Control Lett., 65(2014), 106--111.

\bibitem{Guo-Zhang} B.Z. Guo, L. Zhang, {\it
Local exact controllability to positive trajectory for parabolic system of chemotaxis},
Math. Control Relat. Fields, 6(2016), 143--165.

\bibitem{H2003}   D. Horstmann, {\it
From 1970 until present: the Keller-Segel model in chemotaxis and its consequences. I},
Jahresber. Deutsch. Math.-Verein., 105(2003), 103--165.

\bibitem{HW2005} D. Horstmann,  M. Winkler, {\it
Boundedness vs. blow-up in a chemotaxis system},
J. Differential Equations, 215(2005), 52--107.

\bibitem{JLW2013} H.Y. Jin, J.Y. Li, Z.A. Wang, {\it Asymptotic stability of traveling waves of a chemotaxis model with singular sensitivity},
J. Differential Equations, 255(2013), 193--219.

\bibitem{KS1971} E.F. Keller, L.A. Segel, {\it Traveling bands of chemotactic bacteria: a theoretical analysis,} J. Theoret. Biol., 26(1971), 235--248.

\bibitem{LSU1968} O.A. Lady\v{z}enskaja, V.A. Solonnikov, N.N. Ural'ceva, {\it Linear and quasilinear
equations of parabolic type}, Translations of Mathematical Monographs, Vol. 23 American Mathematical Society, Providence, R.I. 1968.



\bibitem{levine2} H.A. Levine, B.D. Sleeman, {\it
A system of reaction diffusion equations arising in the theory of reinforced random walks},
SIAM J. Appl. Math., 57(1997), 683--730.

\bibitem{levine3} H.A. Levine, B.D. Sleeman, M. Nilsen-Hamilton, {\it
A mathematical model for the roles of pericytes and macrophages in the initiation of angiogenesis. I. The role of protease inhibitors in preventing angiogenesis},
Math. Biosci., 168(2000), 77--115.

\bibitem{LZ2015} H.C. Li, K. Zhao, {\it
Initial-boundary value problems for a system of hyperbolic balance laws arising from chemotaxis},
J. Differential Equations, 258(2015), 302--338.

\bibitem{shi}   P. Lissy, E. Zuazua, {\it
Internal observability for coupled systems of linear partial differential equations},
SIAM J. Control Optim., 57(2019), 832--853.



\bibitem{MN1980}  A. Matsumura, T. Nishida, {\it
The initial value problem for the equations of motion of viscous and heat-conductive gases},
J. Math. Kyoto Univ., 20(1980), 67--104.

\bibitem{2013187} K. Mauffrey, {\it On the null controllability of a $3\times 3$ parabolic system with non--constant coefficients
by one or two control forces}, J. Math. Pures Appl., 99(2013), 187--210.


\bibitem{rwwzz} L.G. Rebholz, D.H. Wang, Z.A. Wang, C. Zerfas, K. Zhao
{\it Initial boundary value problems for a system of parabolic conservation laws arising from chemotaxis in multi-dimensions},
Discrete Contin. Dyn. Syst., 39(2019), 3789--3838.


\bibitem{liu}  D. Steeves, B. Gharesifard, A.-R. Mansouri, {\it
Controllability of coupled parabolic systems with multiple underactuations, Part 2: Null controllability},
SIAM J. Control Optim., 57(2019),  3297--3321.

\bibitem{TY2018}  Q. Tao, Z.A. Yao, {\it
Global existence and large time behavior for a two-dimensional chemotaxis-shallow water system},
J. Differential Equations, 265(2018), 3092--3129.

\bibitem{TWW2013} Y.S. Tao, L.H. Wang,  Z.A. Wang, {\it
Large-time behavior of a parabolic-parabolic chemotaxis model with logarithmic sensitivity in one dimension},
Discrete Contin. Dyn. Syst. Ser. B, 18(2013), 821--845.



\bibitem{T1982}
R. Temam,  {\it Behaviour at time $t= 0$ of the solutions of semilinear evolution equations},
J. Differential Equations, 43(1982), 73--92.



%
%
%
%
%
%
%
%









\bibitem{WXY2016} Z.A. Wang, Z.Y. Xiang, P. Yu, {\it
Asymptotic dynamics on a singular chemotaxis system modeling onset of tumor angiogenesis},
J. Differential Equations, 260(2016), 2225--2258.


\bibitem{W2018} M. Winkler, {\it Renormalized radial large-data solutions to the higher-dimensional Keller-Segel system with singular sensitivity and signal absorption}, J. Differential Equations, 264(2018), 2310--2350.




























\end{thebibliography}
\end{document}